\documentclass[reqno]{amsart}
\allowdisplaybreaks

\usepackage{amssymb}
\usepackage{amsmath,amssymb,amsfonts,amsthm,latexsym,amscd,epsfig,color}
\usepackage{mathrsfs}
\usepackage[font=small]{caption}
\usepackage{array}

\usepackage{a4wide}
\usepackage[T1]{fontenc}
\usepackage[latin1]{inputenc}
\usepackage{tikz}
\usepackage{appendix}
\usetikzlibrary{shapes}

\newtheorem{theo}{Theorem}%[section]
\newtheorem{prop}{Proposition}%[section]
\newtheorem{lemma}{Lemma}%[section]
\newtheorem{cor}{Corollary}
\theoremstyle{definition}
\newtheorem{df}{Definition}%[section]
\newtheorem{ex}{Example}

\theoremstyle{remark}

\newtheorem{rem}{Remark\!\!}

\newtheorem{nrem}{Remark}

\newcommand{\rcp}[1]{\frac{1}{#1}}

\newcommand{\bt}{\begin{theo}}
\newcommand{\et}{\end{theo}}
\newcommand{\bl}{\begin{lemma}}
\newcommand{\el}{\end{lemma}}
\newcommand{\bp}{\begin{prop}}
\newcommand{\ep}{\end{prop}}
\newcommand{\bdf}{\begin{df}}
\newcommand{\edf}{\end{df}}
\newcommand{\brem}{\begin{rem}}
\newcommand{\erem}{\end{rem}}
\newcommand{\bnrem}{\begin{nrem}}
\newcommand{\enrem}{\end{nrem}}
\newcommand{\bex}{\begin{ex}}
\newcommand{\eex}{\end{ex}}
\newcommand{\bcor}{\begin{cor}}
\newcommand{\ecor}{\end{cor}}
\newcommand{\bncor}{\begin{ncor}}
\newcommand{\encor}{\end{ncor}}
\newcommand{\bpf}{\begin{proof}}
\newcommand{\epf}{\end{proof}}

\begin{document}

\title[Exact Counting Results and Corrections]{Counting phylogenetic networks with few reticulation vertices: exact enumeration and corrections}
\author{Michael Fuchs \and Bernhard Gittenberger \and Marefatollah Mansouri}
\thanks{This research has been supported by a bilateral Austrian-Taiwanese project FWF-MOST, grants
I~2309-N35 (FWF) and MOST-104-2923-M-009-006-MY3 (MOST)}
\address{Department of Mathematical Sciences, National Chengchi University,
Taipei, 116, Taiwan.}
\email{mfuchs@nccu.edu.tw}

\address{Department of Discrete Mathematics and Geometry, Technische
Universit\"at Wien, Wiedner Hauptstra\ss e 8-10/104, A-1040 Wien, Austria.}
\email{gittenberger@dmg.tuwien.ac.at}

\address{Department of Discrete Mathematics and Geometry, Technische
Universit\"at Wien, Wiedner Hauptstra\ss e 8-10/104, A-1040 Wien, Austria.}
\email{marefatollah.mansouri@tuwien.ac.at}

%\subjclass{Primary: ???, Secondary: ???}

%\keywords

\date{\today}

\begin{abstract}
In previous work, we gave asymptotic counting results for the number of tree-child and normal networks with $k$ reticulation vertices and explicit exponential generating functions of the counting sequences for $k=1,2,3$. The purpose of this note is two-fold. First, we make some corrections to our previous approach which overcounted the above numbers and thus gives erroneous exponential generating functions (however, the overcounting does not affect our asymptotic counting results). Secondly, we use our (corrected) exponential generating functions to derive explicit formulas for the number of tree-child and normal networks with $k=1,2,3$ reticulation vertices. This re-derives recent results of Carona and Zhang, answers their question for normal networks with $k=2$, and adds new formulas in the case $k=3$.
\end{abstract}

\maketitle

\section{Introduction and Results}
Phylogenetic networks have become a standard tool in evolutionary biology over the last decades since, in contrast to phylogentic trees, they are also able to model reticulation events such as hybridization and lateral gene transfer; see Huson et al. \cite{HRS} and Chapter 10 in Steel \cite{St}.

We will start by defining them. A (binary) {\it phylogenetic network} is a rooted connected DAG
(directed acyclic graph) without double edges whose vertices belong to one of the following four sets:
\begin{itemize}
\item[(a)] A unique root which has indegree $0$ and outdegree $2$;
\item[(b)] Leaves of indegree $1$ and outdegree $0$;
\item[(c)] {\it Tree vertices} of indegree $1$ and outdegree $2$;
\item[(d)] {\it Reticulation vertices} of indegree $2$ and outdegree $1$.
\end{itemize}
Moreover, we say that a phylogenetic network is {\it leaf-labeled} if the leaves are bijectively labeled and {\it vertex-labeled} if all vertices are bijectively labeled.

Next, we recall the two subclasses of phylogenetic networks which were investigated in \cite{FuGiMa}. The first subclass consists of {\it tree-child networks} which are phylogenetic networks where every reticulation vertex is not directly followed by another reticulation vertex and every tree vertex has at least one child which is not a reticulation vertex; see Cardona et al. \cite{GFG}. The second subclass are {\it normal networks} which are tree-child networks with the additional constraint that if there is a (directed) path between two vertices of length at least $2$, then there is no direct edge between these vertices; see Willson \cite{Wi1,Wi2}.

Finally, we recall the following notations:
\begin{center}
\begin{tabular}{ll}
$\tilde{T}_{k,\ell}$ & number of leaf-labeled tree-child networks with $\ell$ leaves and $k$ reticulation vertices;\\
$T_{k,n}$ & number of vertex-labeled tree-child networks with $n$ vertices and $k$ reticulation vertices;\\
$\tilde{N}_{k,\ell}$ & number of leaf-labeled normal networks with $\ell$ leaves and $k$ reticulation vertices;\\
$N_{k,n}$ & number of vertex-labeled normal networks with $n$ vertices and $k$ reticulation vertices.
\end{tabular}
\end{center}

The main purpose of \cite{FuGiMa} was the proof of the following asymptotic counting result: for fixed $k$, we have
\begin{equation}\label{asym-vertex}
T_{k,n}\sim N_{k,n}\sim c_k\left(1-(-1)^n\right)\left(\frac{\sqrt{2}}{e}\right)^nn^{n+2k-1},\qquad (n\rightarrow\infty),
\end{equation}
where $c_k$ is a computable constant. Moreover, also in \cite{FuGiMa}, a similar result for the leaf-labeled cases was derived from this as a consequence: for fixed $k$, we have
\begin{equation}\label{asym-leaf}
\tilde{T}_{k,\ell}\sim \tilde{N}_{k,\ell}\sim 2^{3k-1}c_k\left(\frac{2}{e}\right)^{\ell}{\ell}^{\ell+2k-1},\qquad (\ell\rightarrow\infty).
\end{equation}

In the recent paper \cite{Zh2}, Cardona and Zhang introduced an algorithmic method which can be used efficiently to compute the values of $\tilde{T}_{k,\ell}$ if $k$ and $\ell$ are small. Moreover, they showed that their method also yields formulas for $\tilde{T}_{k,\ell}$ for $k=1,2$ and all values of $\ell$. (The formula for $k=1$ was also contained in Zhang \cite{Zh1}.)

One purpose of this note, is to point out that the formulas for $\tilde{T}_{k,\ell}$ with $k=1, k=2$
and even $k=3$ also follow from our results in \cite{FuGiMa}. This is, because for these values of $k$ we gave explicit expressions for the exponential generating function $T_k(z)$ of $T_{k,n}$ in \cite{FuGiMa}. More precisely, we showed in \cite{FuGiMa} that for fixed $k$:
\begin{equation}\label{egf-T}
T_{k}(z)=z\frac{\tilde{a}_k^{[T]}(z^2)-\tilde{b}_k^{[T]}(z^2)\sqrt{1-2z^2}}{(1-2z^2)^{2k-1/2}}
\end{equation}
with polynomials $\tilde{a}_k^{[T]}(z)$ and $\tilde{b}_k^{[T]}(z)$ which we computed in \cite{FuGiMa} for $k=1,2,3$; see below for corrections for the expressions from \cite{FuGiMa}. (In principle, our method can also be used to compute these polynomials for higher values of $k$ but the computation becomes more and more cumbersome.) From this, we obtain
\begin{align}
T_{k,2n+1}&=(2n+1)![z^{2n+1}]T_{k}(z)\nonumber\\
&=(2n+1)![z^n]\frac{\tilde{a}_k^{[T]}(z)-\tilde{b}_k^{[T]}(z)\sqrt{1-2z}}{(1-2z)^{2k-1/2}}\nonumber\\
&=(2n+1)!\left(2^{-n}\binom{2n}{n}\tilde{r}_k^{[T]}(n)-2^{n}\tilde{p}_k^{[T]}(n)\right),\label{exp-T}
\end{align}
where $\tilde{r}_k^{[T]}(n)$ is a rational function in $n$ and $\tilde{p}_k^{[T]}(n)$ is a polynomial in $n$; see the next section for details and explicit expressions for $\tilde{r}_k^{[T]}(n)$ and $\tilde{p}_k^{[T]}(n)$ when $k=1,2,3$. Then, from the equation (see \cite{FuGiMa})
\[
\tilde{T}_{k,\ell}=\frac{\ell!}{(2\ell+2k-1)!}T_{k,2\ell+2k-1},
\]
we also obtain explicit results for $\tilde{T}_{k,\ell}$ when $k=1,2,3$. Our formula for $k=2$ slightly simplifies the one given in \cite{Zh2} and the formula for $k=3$ correctly produces all the terms given for $\tilde{T}_{3,\ell}$ in Table~3 in \cite{Zh2}.

Similarly, we obtain explicit expressions for $N_{k,2n+1}$ and $\tilde{N}_{k,\ell}$ for $k=1,2,3$ since for these cases, we again have explicit results for the exponential generating function $N_{k}(z)$ of $N_{k,n}$ which has the same form as that of tree-child networks:
\begin{equation}\label{egf-N}
N_{k}(z)=z\frac{\tilde{a}_k^{[N]}(z^2)-\tilde{b}_k^{[N]}(z^2)\sqrt{1-2z^2}}{(1-2z^2)^{2k-1/2}},
\end{equation}
where $\tilde{a}_k^{[N]}(z)$ and $\tilde{b}_k^{[N]}(z)$ are polynomials which where derived for $k=1,2,3$ in \cite{FuGiMa} (again see below for corrections). Thus, as above,
\begin{equation}\label{exp-N}
N_{k,2n+1}=(2n+1)!\left(2^{-n}\binom{2n}{n}\tilde{r}_k^{[N]}(n)-2^{n}\tilde{p}_k^{[N]}(n)\right)
\end{equation}
with a rational function $\tilde{r}_k^{[N]}(n)$ and a polynomial $\tilde{p}_k^{[N]}(n)$ which are again given for $k=1,2,3$ in the next section. Moreover,
\[
\tilde{N}_{k,\ell}=\frac{\ell!}{(2\ell+2k-1)!}N_{k,2\ell+2k-1}.
\]
The explicit formula for $\tilde{N}_{1,\ell}$ was already given in \cite{Zh1} and finding an
explicit formula for $\tilde{N}_{2,\ell}$ was posed as an open problem in \cite{Zh2}. Our formula
for $\tilde{N}_{2,\ell}$ correctly produces all the corresponding values from Table~1 in \cite{Zh2}. Moreover, our formula for $\tilde{N}_{3,\ell}$ correctly produces the first two values in that table and corrects the remaining two. (See the end of Section~\ref{normal} for details.)

When using the above mentioned results from \cite{FuGiMa} to derive the above formulas, our values
for $k=2$ and $k=3$ initially differed from those given in \cite{Zh2,Zh1}. The reason for this is
that we forgot to consider some cases in \cite{FuGiMa}. (These cases are asymptotically negligible
and thus do not affect the main results from \cite{FuGiMa} displayed in \eqref{asym-vertex} and
\eqref{asym-leaf}; however, they do affect the exponential generating functions for $k=2,3$ from
\cite{FuGiMa}.) Thus, the second purpose of this note is to explain what we forgot and give the
correct expressions for the polynomials $\tilde{a}_k^{[\star]}(z)$ and $\tilde{b}_k^{[\star]}(z)$
with $k=2,3$ and $\star\in\{T,N\}$. We collect them in the next two theorems, where, for the sake of completeness, we also include the expressions for $k=1$.

\bt[Tree-Child Networks]\label{polys-T}
The polynomials $\tilde{a}_k^{[T]}(z)$ and $\tilde{b}_k^{[T]}(z)$ in the expression (\ref{egf-T}) for $T_k(z)$ with $k=1,2,3$ are as follows.
\begin{itemize}
\item[(i)] $\tilde{a}_1^{[T]}(z)=z$ and $\tilde{b}_1^{[T]}(z)=z$;
\item[(ii)] $\tilde{a}_2^{[T]}(z)=-z^4+8z^3$ and $\tilde{b}_2^{[T]}(z)=8z^3$;
\item[(iii)] $\tilde{a}_3^{[T]}(z)=-35z^6+175z^5$ and $\tilde{b}_3^{[T]}(z)=34z^6+175z^5$.
\end{itemize}
\et

\bt[Normal Networks]\label{polys-N}
The polynomials $\tilde{a}_k^{[N]}(z)$ and $\tilde{b}_k^{[N]}(z)$ in the expression (\ref{egf-N}) for $N_k(z)$ with $k=1,2,3$ are as follows.
\begin{itemize}
\item[(i)] $\tilde{a}_1^{[N]}(z)=2-3z$ and $\tilde{b}_1^{[N]}(z)=2-z$;
\item[(ii)] $\tilde{a}_2^{[N]}(z)=11z^4-66z^3+50z^2-8z$ and $\tilde{b}_2^{[N]}(z)=-28z^3+42z^2-8z$;
\item[(iii)] $\tilde{a}_3^{[N]}(z)=877z^6-3065z^5+2392 z^4-628 z^3+64z^2 $ and $\tilde{b}_3^{[N]}(z)=110z^6-1455z^5+1860z^4-564 z^3+64z^2 $.
\end{itemize}
\et

We conclude the introduction by a short sketch of this note. In the next section, we will give
more details on the derivation of \eqref{exp-T} and \eqref{exp-N} and list the expressions for $\tilde{r}_k^{[\star]}(n)$ and $\tilde{p}_k^{[\star]}(n)$ for $k=1,2,3$ and $\star\in\{T,N\}$. In Section~\ref{meth-FuGiMa}, we will recall the method from \cite{FuGiMa}. Then, we will explain in
Sections~\ref{tree-child} and~\ref{normal} how the method is corrected to yield the results from
Theorem~\ref{polys-T} for tree-child networks and Theorem~\ref{polys-N} for normal networks,
respectively. Finally, an appendix will contain the answer of a counting problem which is of independent interest and can be used to have a quick verification of whether the values produced by the above formulas (and values published in other works) are reasonable or not.
\begin{table}[ht]
\centering
\begin{tabular}{|l|m{1.7cm}|m{1.6cm}|m{2.2cm}|}
\hline
\qquad\qquad\qquad\qquad Type of networks &  EGFs in \cite{FuGiMa} &Corrected EGFs& Formulas \\
\hline
Tree-child networks with one reticulation vertex & Prop. $4.1$ &  & \cite{Zh2,Zh1}, Thm.~\ref{Ex1} \\  \cline{2-2}\cline{4-4}
Tree-child networks with two reticulation vertices & Prop. $4.2$ & Thm.~\ref{polys-T} & \cite{Zh2}, Thm.~\ref{Ex1}   \\ \cline{2-2}\cline{4-4}
Tree-child networks with three reticulation vertices& Prop. $4.3$ & & Thm.~\ref{Ex1}\\
\hline \hline
Normal networks with one reticulation vertex& Prop. $3.1$ & & \cite{Zh1}, Thm.~\ref{Ex2} \\\cline{2-2}\cline{4-4}
Normal networks with two reticulation vertices & Prop. $3.2$ & Thm.~\ref{polys-N} & Thm.~\ref{Ex2} \\ \cline{2-2}\cline{4-4}
Normal networks with three reticulation vertices & Prop. $3.3$ & & Thm.~\ref{Ex2}\\
\hline
\end{tabular}
\smallskip
\caption{Overview of the main results. EGF means exponential generating function. The results in the second column appeared in \cite{FuGiMa}. The third column contains
the corrected results from the current paper. (The results in the cases of one reticulation vertex have already been correct in \cite{FuGiMa}.) In the fourth column, previous known explicit formulas are listed together with the formulas derived in this paper. The formula for normal networks with two reticulation vertices was mentioned in \cite{Zh2} as an open problem.
\label{tab:overview}}
  \end{table}

\section{Explicit Formulas for the Number of Tree-Child and Normal Networks with $k=1,2,3$.}

In this section, we fill in the missing steps for (\ref{exp-T}) and (\ref{exp-N}). More precisely, we give more details for the last equality in (\ref{exp-T}). Therefore, we drop the superscript and thus consider
\[
[z^n]\frac{\tilde{a}_k(z)-\tilde{b}_k(z)\sqrt{1-2z}}{(1-2z)^{2k-1/2}}.
\]
First, note that
\[
[z^n]\frac{z^m}{(1-2z)^{\alpha}}=2^{n-m}\binom{-\alpha}{n-m}(-1)^{n-m}=2^{n-m}\binom{n-m+\alpha-1}{n-m}.
\]
Using this gives
\begin{align*}
[z^n]&\frac{\tilde{a}_k(z)-\tilde{b}_k(z)\sqrt{1-2z}}{(1-2z)^{2k-1/2}}\\
&=\sum_{m\geq 0}([z^m]\tilde{a}_k(z))2^{n-m}\binom{n-m+2k-3/2}{n-m}-\sum_{m\geq 0}([z^m]\tilde{b}_k(z))2^{n-m}\binom{n-m+2k-2}{n-m}.
\end{align*}

For the second term, we have
\[
\binom{n-m+2k-2}{n-m}=\binom{n-m+2k-2}{2k-2}
\]
which is a polynomial in $n$ and thus
\[
\tilde{p}_k(n)=\sum_{m\geq 0}([z^m]\tilde{b}_k(z))2^{-m}\binom{n-m+2k-2}{2k-2}
\]
is also a polynomial in $n$.

For the first term, observe that
\[
\binom{n-m+2k-3/2}{n-m}=4^{-n}\binom{2n}{n}\tilde{r}_{k,m}(n)
\]
with a suitable rational function $\tilde{r}_{k,m}(n)$ in $n$ whose coefficients depend on $k$ and $m$. Thus,
\[
\tilde{r}_k(n)=\sum_{m\geq 0}([z^m]\tilde{a}_k(z))2^{-m}\tilde{r}_{k,m}(n)
\]
is also a rational function in $n$.

Collecting everything gives now
\[
[z^n]\frac{\tilde{a}_k(z)-\tilde{b}_k(z)\sqrt{1-2z}}{(1-2z)^{2k-1/2}}=2^{-n}\binom{2n}{n}\tilde{r}_k(n)-2^{n}\tilde{p}_k(n)
\]
which is the claimed form in (\ref{exp-T}) and (\ref{exp-N}).

From Theorem~\ref{polys-T} and Theorem~\ref{polys-N} and some computation, we now can find explicit expressions for $\tilde{r}_k^{[\star]}(n)$ and $\tilde{p}_k^{[\star]}(n)$ for $k=1,2,3$ and $\star\in\{T,N\}$ and thus have explicit formulas for $T_{k,2n+1},\tilde{T}_{k,\ell},N_{k,2n+1}$ and $\tilde{N}_{k,\ell}$ for $k=1,2,3$.

\bt[Tree-Child Networks] \label{Ex1}
The rational function $\tilde{r}_k^{[T]}(n)$ and the polynomial $\tilde{p}_k^{[T]}(n)$ in the
formula (\ref{exp-T}) for $k=1,2,3$ are as follows:
\begin{itemize}
\item[(i)] $\tilde{r}_1^{[T]}(n)=n$ and $\tilde{p}_1^{[T]}(n)=\rcp{2}$;
\item[(ii)] $\tilde{r}_2^{[T]}(n)=\dfrac{n(n-1)(n-2)(3n-1)}{3(2n-1)}$ and $\tilde{p}_2^{[T]}(n)=\rcp{2}(n-1)(n-2)$;
\item[(iii)] $\tilde{r}_3^{[T]}(n)= \dfrac{n^2(n-1)(n-2)(n-3)(n-4)}{3(2n-1)}$ and $\tilde{p}_3^{[T]}(n)=\frac{1}{192}(n-2)(n-3)(n-4)(48n-65)$.
\end{itemize}
\et

\bt[Normal Networks]\label{Ex2}
The rational function $\tilde{r}_k^{[N]}(n)$ and the polynomial $\tilde{p}_k^{[N]}(n)$ in the
formula (\ref{exp-N}) for $k=1,2,3$ are as follows:
\begin{itemize}
\item[(i)] $\tilde{r}_1^{[N]}(n)=(n+2)$ and $\tilde{p}_1^{[N]}(n)=\frac{3}{2}$;
\item[(ii)] $\tilde{r}_2^{[N]}(n)=\dfrac{n(3n-7)(n^2+9n-4)}{3(2n-1)}$ and $\tilde{p}_2^{[N]}(n)=\dfrac{(n+1)(3n-7)}{2}$;
\item[(iii)] $\tilde{r}_3^{[N]}(n)=\dfrac{n(n-1)}{3(2n-1)}(n^4+15n^3-158n^2+324n+40)$ and $\tilde{p}_3^{[N]}(n)=\frac{1}{192}(144n^4-751n^3-1089n^2+9106n-7080)$.
\end{itemize}
\et

\section{Summary of the Method from \cite{FuGiMa}}\label{meth-FuGiMa}

In this section, we recall the method from \cite{FuGiMa}.

\begin{center}
  	\begin{figure}[h]
  		\begin{minipage}{1\textwidth}
  			\begin{center}
  				{\includegraphics[width=.8\textwidth]{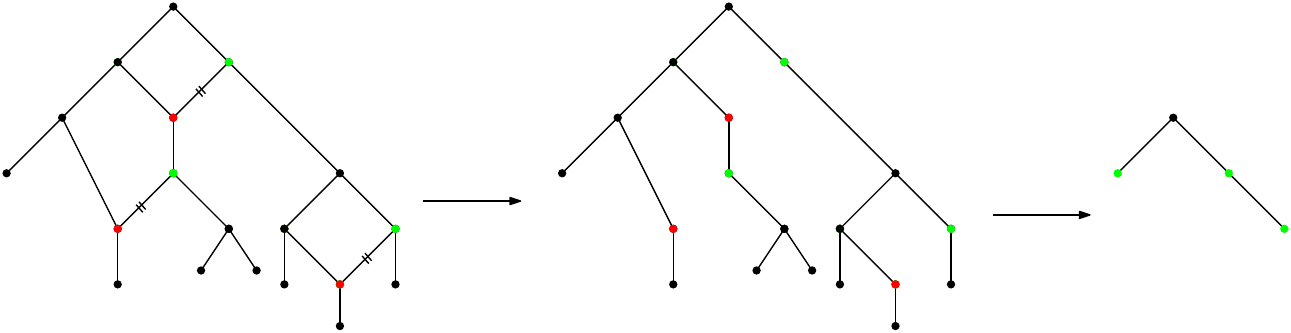}}
  			\end{center}
  		\end{minipage}
  		\caption{A tree-child network (which is even normal) together with a possible choice of a colored Motzkin skeleton and coresponding
  			sparsened skeleton.}	
  		\label{redgreen}
  	\end{figure}
  \end{center}

First, fix a vertex-labeled tree-child network with $k$ reticulation vertices; see Figure~\ref{redgreen} for an example where we dropped all labels and directions are from the root downward. Then, in \cite{FuGiMa}, we performed the following two steps.
\begin{enumerate}
\item[(i)] Color all reticulation vertices red and for each reticulation vertex, pick an incoming edge and color its parent green; then remove the picked edges. Note that the resulting graph is a Motzkin tree (i.e., a rooted and vertex-labeled tree with binary vertices, unary vertices and leaves) with exactly $k$ red and $k$ green unary vertices. We called in \cite{FuGiMa} this Motzkin tree a {\it colored Motzkin skeleton} of the given tree-child network; see Figure~\ref{redgreen} for an example.
\item[(ii)] Contract paths (and the subtrees dangling from them) between green vertices and between green vertices and their last common ancestor, and also remove trees below green vertices in the colored Motzkin skeleton so that a new Motzkin tree is obtained which describes the ancestral relationship of the green vertices; see again Figure~\ref{redgreen} for an example. We called in \cite{FuGiMa} this new Motzkin tree the {\it sparsened skeleton} of the colored Motzkin skeleton.
\end{enumerate}
Now for the construction of all vertex-labeled tree-child or normal networks with $k$ reticulation vertices, we reversed the above process. More precisely, we first considered all possible sparsened skeletons. Picking one of them, we added back the removed paths from step (ii) above and the trees below the green vertices. Here, we worked with the symbolic method from \cite{FlSe} and multivariate exponential generating functions.

More precisely, first consider Motzkin trees which satisfy the tree-child condition with unary vertices playing the role of the reticulation vertices. Such trees are counted with the exponential generating function $M(z,y)$ where $y$ marks unary vertices and $z$ marks all vertices ($M(y,z)$ is only exponential in $z$ and ordinary in $y$):
\[
M(z,y)=\frac{(1+yz)\left(1-\sqrt{1-2z^2-4yz^3}\right)}{z(1+2yz)};
\]
see \cite{FuGiMa} for details. Then, this exponential generating function is used to construct the above mentioned contracted paths with subtrees dangling from them where again the tree-child condition for unary vertices must hold. For this, in \cite{FuGiMa}, we used for tree-child networks the exponential generating function $\hat{P}(z,y,\tilde{y},\hat{y})$ (again exponential in $z$ and ordinary in $y$):
\[
\hat{P}(z,y,\tilde{y},\hat{y})=\frac{1+z\hat{y}}{1-zM(z,\tilde{y})-z^2y\tilde{M}(z,\tilde{y})},
\]
where $y,\hat{y}$ mark unary vertices on the path with $\hat{y}$ the first vertex on the path and $\tilde{y}$ marks unary vertices in the subtrees dangling from the path. (Here, $\tilde{M}(z,y)$ counts those labeled Motzkin trees from $M(z,y)$ which do not start with a unary vertex.) For normal networks, the above exponential generating function had to be replaced by
\begin{equation}\label{path-nn}
P(z,y,\tilde{y},\hat{y})=\frac{1+z\hat y}{1-(z+2z^2y)\tilde M(z,\tilde y)},
\end{equation}
where $\tilde{M}(z,y)$ and $\hat{y}$ are as above, and $y$ now marks unary vertices on the path or which are children of vertices on the path and $\tilde{y}$ marks the remaining unary vertices; see \cite{FuGiMa}.

Now, in \cite{FuGiMa}, we used the above exponential generating functions to construct all colored Motzkin skeletons for tree-child networks resp. normal networks. Finally, we added back the edges from the green vertices to the red vertices by pointing (which on the level of generating functions corresponds to differentiation).
  \begin{center}
  	\begin{figure}[h]
  		\begin{minipage}{.8\textwidth}
  			\begin{center}
  				{\includegraphics[width=1\textwidth]{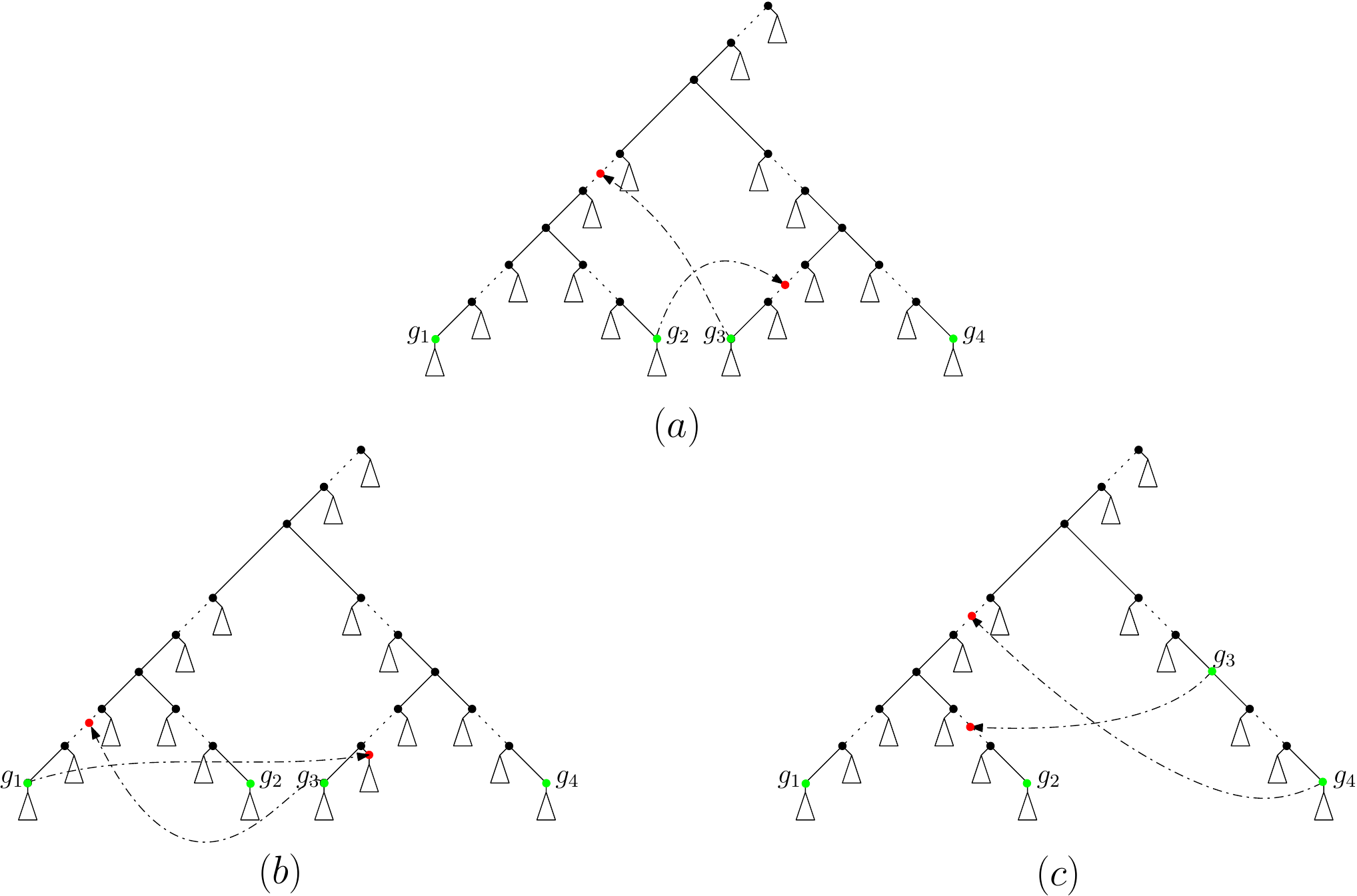}}
  			\end{center}
  		\end{minipage}
  		\caption{ $(a)$: The pointing of $g_2$ and $g_3$ creates the type of cycle which caused overcounting in \cite{FuGiMa}; $(b)$ and $(c)$: The two possible types of near-cycles which caused overcounting in \cite{FuGiMa} when counting normal networks.}	
  		\label{cycle}
  	\end{figure}
  \end{center}
For tree-child networks, the last step seems to be easy, because the above constructions already guarantee that the tree-child condition will hold. However, one has to be careful not to create cycles which could happen if a green vertex points on a vertex on a path from the root leading to it or (less obviously) if for two green vertices $g_i,g_j$ which are not on a same path, $g_i$ points on a vertex on a path leading from the last common ancestor of $g_i$ and $g_j$ to $g_j$, and $g_j$ points on a vertex on a path leading from the last common ancestor of $g_i$ and $g_j$ to $g_i$; see Figure~\ref{cycle}-(a) for an example with $i=2$ and $j=3$. We forgot to subtract networks containing the second type of cycles in \cite{FuGiMa} and thus the generating functions $T_k(z)$ in \cite{FuGiMa} overcounts the number of tree-child networks with $k$ reticulation vertices (however, asymptotically this overcount is not relevant). We will show in the next section how to modify our approach from \cite{FuGiMa} to avoid counting these additional networks containing these cycles.

For normal networks, even more care has to be taken in the above pointing step since, in addition
to cycles, one also has now to be careful not to create {\it near-cycles} by which we mean closed paths
where all edges except one are in the same direction (this is forbidden by the definition of normal
networks). In \cite{FuGiMa} the creation of most of these near-cycles was avoided, however, we
missed the following two more subtle ones: (i) if for two green vertices $g_i,g_j$ which are not
on the same path, $g_i$ points on a vertex on a path leading from the last common ancestor of
$g_i$ and $g_j$ to $g_j$, and $g_j$ points on a child of a vertex on a path leading from the
last common ancestor of $g_i$ and $g_j$ to $g_i$ (or vice versa) and (ii) if for two green
vertices $g_i,g_j$ which are on the same path, both vertices point at vertices on a path and the
pointers cross each other; see Figure~\ref{cycle}-(b) ($i=2$ and $j=1$) and (c) ($i=3$ and $j=4$)
for a depiction of these two cases. We will explain in Section~\ref{normal} how to modify our
approach from \cite{FuGiMa} to avoid counting networks with these kinds of near-cycles as well as
the type of cycles from Figure~\ref{cycle}-(a), which we also forgot to rule out for normal networks in \cite{FuGiMa}.

\section{Tree-Child Networks with $2$ and $3$ Reticulation Vertices}\label{tree-child}

Here, we will give details for the counting of tree-child networks. Since, as explained in the previous section, we overcounted them in \cite{FuGiMa}, we could just take our results from \cite{FuGiMa} and subtract the networks containing the kind of cycles described in the last section. Alternatively, we can start from the scratch and count these networks so that the occurrence of these cycles is avoided in the first place. We will explain the second approach here (where subtractions are, however, still necessary in some cases).

\medskip
\subsection{Tree-Child Networks with Two Reticulation Vertices}
\begin{center}
  	\begin{figure}[h]
  		\begin{minipage}{1\textwidth}
  			\begin{center}
  				{\includegraphics[width=0.4\textwidth]{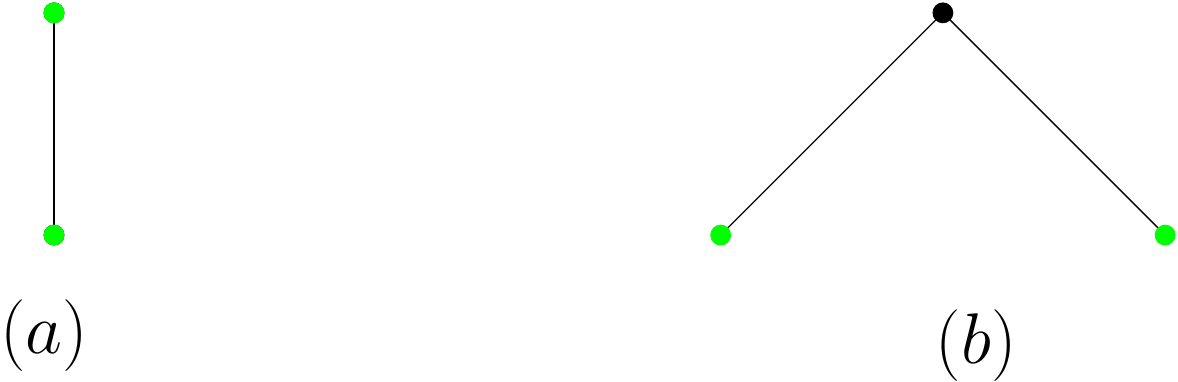}}
  			\end{center}
  		\end{minipage}
  		\caption{All possible sparsened skeletons with two reticulation vertices.% (green/light-gray vertices).
  			}
  		\label{spsk-k2}
  	\end{figure}
\end{center}

 All possible sparsened skeletons are listed in Figure~\ref{spsk-k2}. Note that the creation of the type of cycles explained in the last section in the pointing step is not possible for the one in Figure~\ref{spsk-k2}-(a). Thus, no overcounting has occurred in \cite{FuGiMa} for that sparsened skeleton and we can thus concentrate on the one in Figure~\ref{spsk-k2}-(b).

\begin{center}
	\begin{figure}
		\begin{minipage}{1.3\textwidth}
			\begin{center}
				{\includegraphics[width=.75\textwidth]{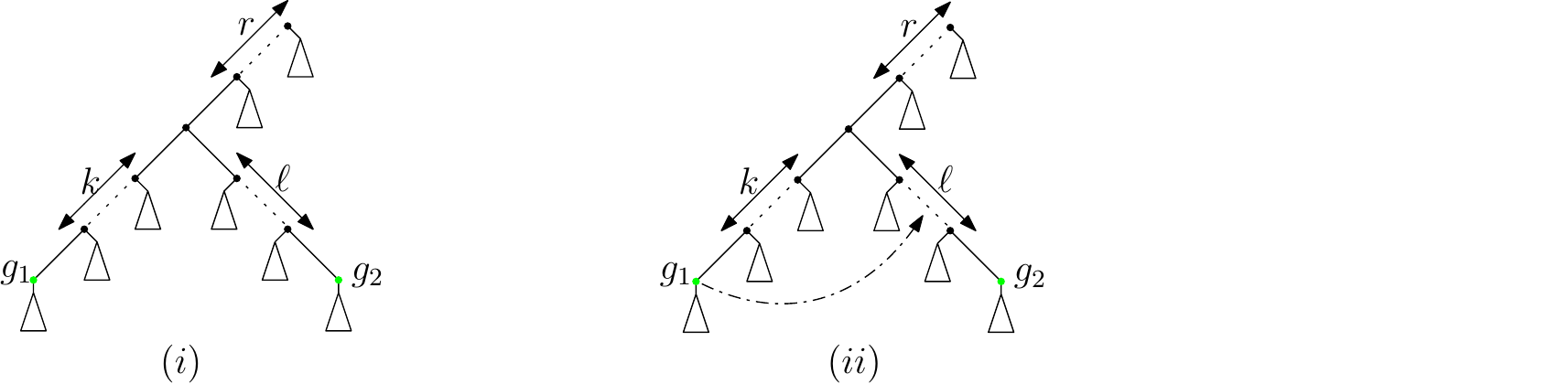}}
			\end{center}
		\end{minipage}
		\caption{The colored Motzkin skeletons arising from the sparsened skeleton in Figure~\ref{spsk-k2}-(b) classified according to the indicated pointing rule; the remaining pointers from green vertices have no restrictions except that they are not allowed to point at a vertex on a path.}
		\label{T2-2}
	\end{figure}
\end{center}

For that one, we have to consider two cases which are listed in Figure~\ref{T2-2}, where the
missing pointers are not allowed to point at vertices on the paths which are contracted (and dangling subtrees deleted) in step (ii) of the construction at the beginning of Section~\ref{meth-FuGiMa}. We explain now in detail the exponential generating functions for these two cases. First, for the networks arising from the colored Motzkin skeletons from Figure~\ref{T2-2}-(i), we have
\[
\frac{1}{2}\partial_{y_1}\partial_{y_2}\left(z^3\tilde{M}(z,y_1+y_2)^2\hat{P}(z,0,y_1+y_2,0)^3\right)\Big\vert_{y_1=y_2=0},
\]
where $y_1$ resp. $y_2$ track possible targets of the pointers starting from $g_1$ resp. $g_2$, the factor $1/2$ comes from symmetry, the factor $z^3$ counts the two green vertices and their last common ancestor, $\tilde{M}(z,y_1+y_2)^2$ counts the two subtrees dangling from the
green vertices (note that pointing at the roots of these subtrees is not allowed and we thus have to use $\tilde{M}$ instead of $M$), and $\hat{P}(z,0,y_1+y_2,0)^3$ counts the three paths $k,\ell,$ and $r$. Next, for the networks arising from the colored Motzkin skeletons from Figure~\ref{T2-2}-(ii), we have
\[
\partial_{y_1}\partial_{y_2}\left(z^3\tilde{M}(z,y_2)^2\hat{P}(z,y_1,y_2,y_1)\hat{P}(z,0,y_2,0)^2\right)\Big\vert_{y_1=y_2=0},
\]
where the terms are explained as above with $\hat{P}(z,y_1,y_2,y_1)$ corresponding to path $\ell$
and $\hat{P}(z,0,y_2,0)^2$ corresponding to the remaining two paths (note that since only $g_1$ is
allowed to point at a vertex on a path, we do not have to consider symmetry).

Now, summing the above two exponential generating functions gives the exponential generating function counting the tree-child networks arising from the sparsened skeleton in Figure~\ref{spsk-k2}-(b). Adding with the exponential generating function of the networks arising from Figure~\ref{spsk-k2}-(a) and dividing the result by $2^k=4$ (since every tree-child network is obtained from this procedure exactly $4$ times), we obtain
\[
T_{2}(z)=z\frac{-z^8+8z^6-8z^6\sqrt{1-2z^2}}{(1-2z^2)^{7/2}}.
\]
This implies
\[
T_{2,2n+1}=(2n+1)!\left(\frac{n(n-1)(n-2)(3n-1)}{3(2n-1)2^{n}}\binom{2n}{n}-2^{n-1}(n-1)(n-2)\right)
\]
and
\begin{equation}
\tilde{T}_{2,\ell}=\ell!\left(\frac{(\ell+1)\ell(\ell-1)(3\ell+2)}{6(2\ell+1)2^{\ell}}\binom{2\ell+2}{\ell+1}-2^{\ell}\ell(\ell-1)\right).
\label{eq:TL2}
\end{equation}
The latter sequence starts with (for $\ell\geq 3$)
\[
42, 1272, 30300,  696600, 16418430,  405755280, \ldots
\]
which matches with the values from Table 3 in \cite{Zh2}.

\medskip
\subsection{Tree-Child Networks with Three Reticulation Vertices}
\begin{center}
  	\begin{figure}[h]
  		\begin{minipage}{1\textwidth}
  			\begin{center}
  				{\includegraphics[width=0.4\textwidth]{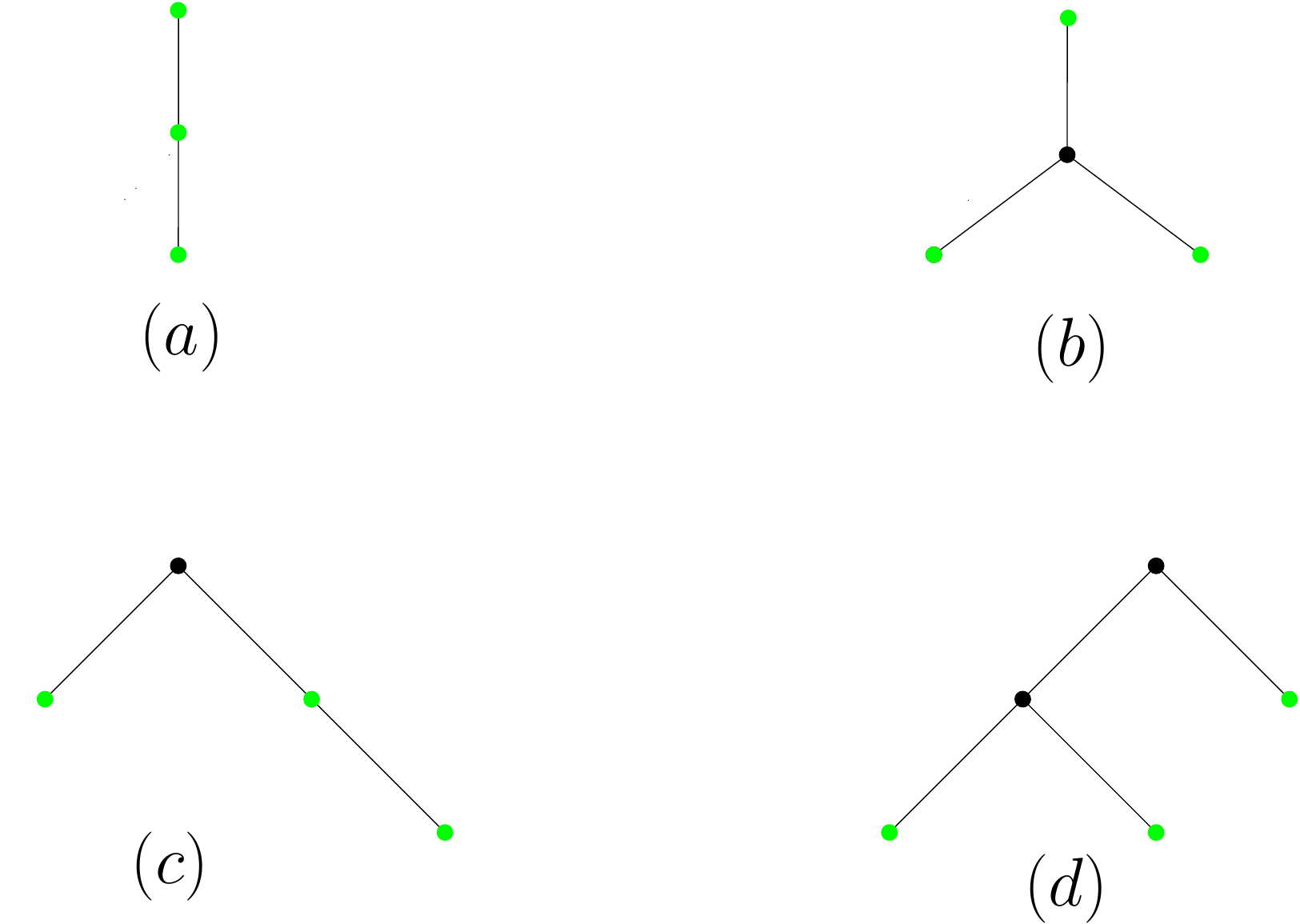}}
  			\end{center}
  		\end{minipage}
  		\caption{All possible sparsened skeletons with three reticulation vertices.% (green/light-gray vertices).
  			}
  		\label{spsk-k3}
  	\end{figure}
\end{center}

Next, we consider three reticulation vertices. All sparsened skeletons
are listed in Figure~\ref{spsk-k3}. As in the case $k=2$, we do not have to consider the sparsened
skeleton in Figure~\ref{spsk-k3}-(a), because the exponential generating function for counting all
tree-child networks arising from it in \cite{FuGiMa} is already correct. All other cases have to
be re-considered, which we will do now.

\begin{center}
  	\begin{figure}[h]
  		\begin{minipage}{1\textwidth}
  			\begin{center}
  				{\includegraphics[width=0.85\textwidth]{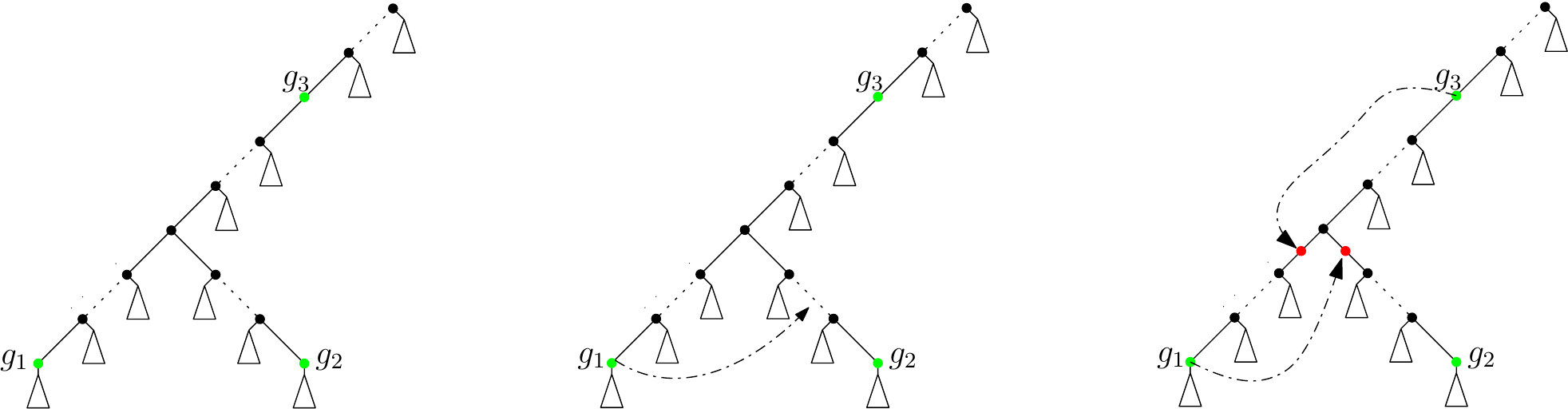}}
  			\end{center}
  		\end{minipage}
  		\caption{The colored Motzkin skeletons arising from the sparsened skeleton in
Figure~\ref{spsk-k3}-(b) classified according to the indicated pointing rules; the missing pointers have no restrictions except that the pointers from $g_1$ and $g_2$ are not allowed to point at vertices on a path. The exponential generating function of the first two colored Motzkin skeletons have to be added, whereas the one from the last one has to be subtracted.}
  		\label{T3-2}
  	\end{figure}
\end{center}

First, for the colored Motzkin skeletons arising from the sparsened skeleton in
Figure~\ref{spsk-k3}-(b), we classify them according to the cases in Figure~\ref{T3-2}. Here, the exponential
generating functions of the first two cases have to be added up, whereas the exponential
generating function of the last case must be subtracted, because $g_1$ and $g_3$ are not allowed
to point to the children of the last common ancestor of $g_1$ and $g_2$, because that ancestor is a tree vertex and thus cannot have two reticulation vertices as children. Overall, we get for the exponential generating function in this case

\begin{align*}
\rcp2{\bf Y}&\left(z^4\tilde{M}(z,y_1+y_2+y_3)^2\hat{P}(z,y_3,y_1+y_2+y_3,y_3)^2\hat{P}(z,y_3,y_1+y_2+y_3,0)\right.\\
&\qquad\quad\left.\times\hat{P}(z,0,y_1+y_2+y_3,0)\right)\\
+&{\bf Y}\left(z^4\tilde{M}(z,y_2+y_3)^2\hat{P}(z,y_1+y_3,y_2+y_3,y_1+y_3)\hat{P}(z,y_3,y_2+y_3,y_3)\right.\\
&\qquad\quad\left.\times\hat{P}(z,y_3,y_2+y_3,0)\hat{P}(z,0,y_2+y_3,0)\right)\\
-&{\bf Y}\left(z^4\tilde M(z,y_2)^2\hat{P}(z,0,y_2,y_1)\hat{P}(z,0,y_2,y_3)\hat{P}(z,0,y_2,0)^2\right),
\end{align*}
where ${\bf Y}(\cdot)$ is used as an abbreviation for $\partial_{y_1}\partial_{y_2}\partial_{y_3}(\cdot)\Big\vert_{y_1=y_2=y_3=0}$.

\begin{center}
 	\begin{figure}[h]
 		\begin{minipage}{1\textwidth}
 			\begin{center}
 				{\includegraphics[width=0.9\textwidth]{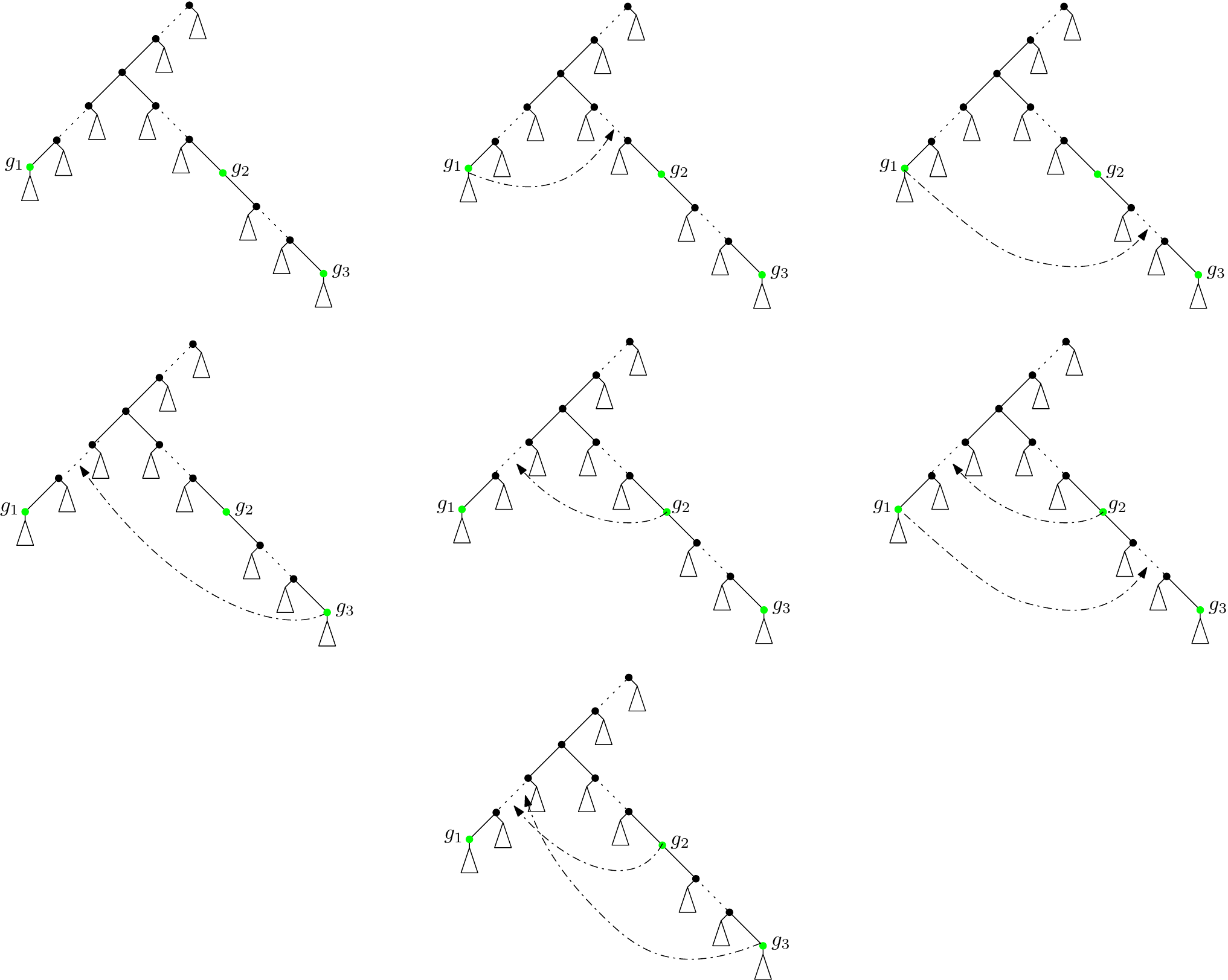}}
 			\end{center}
 		\end{minipage}
 		\caption{The colored Motzkin skeletons arising from the sparsened skeleton in
Figure~\ref{spsk-k3}-(c) classified according to the indicated pointing rules; the missing pointers of green vertices have no restrictions except that $g_1$ and $g_3$ are not allowed to point at vertices on a path and $g_2$ is not allowed to points at a vertex of the path above it.}
 		\label{T3-3}
 	\end{figure}
 \end{center}

Next, we consider the sparsened skeleton in Figure~\ref{spsk-k3}-(c) whose colored Motzkin
skeletons are classified into the cases given in Figure~\ref{T3-3}. From them, we obtain the following equation for the exponential
generating function:
\begin{align*}
{\bf Y}&\left(z^4\tilde M(z,y_1+y_2+y_3)^2\hat{P}(z,0,y_1+y_2+y_3,0)^4\right)\\
+&{\bf Y}\left(z^4\tilde M(z,y_2+y_3)^2\hat{P}(z,y_1,y_2+y_3,y_1)\hat{P}(z,y_2,y_2+y_3,0)\hat{P}(z,0,y_2+y_3,0)^2\right)\\
+&{\bf Y}\left(z^4\tilde M(z,y_2+y_3)^2\hat{P}(z,y_1+y_2,y_2+y_3,y_1)\hat{P}(z,0,y_2+y_3,0)^3\right)\\
+&{\bf Y}\left(z^4\tilde M(z,y_1+y_2)^2\hat{P}(z,y_3,y_1+y_2,y_3)\hat{P}(z,y_2,y_1+y_2,0)\hat{P}(z,0,y_1+y_2,0)^3\right)\\
+&{\bf Y}\left(z^4\tilde M(z,y_1+y_3)^2\hat{P}(z,y_2,y_1+y_3,y_2)\hat{P}(z,0,y_1+y_3,0)^3\right)\\
+&{\bf Y}\left(z^4\tilde M(z,y_3)^2\hat{P}(z,y_2,y_3,y_2)\hat{P}(z,y_1,y_3,0)\hat{P}(z,0,y_3,0)^2\right)\\
+&{\bf Y}\left(z^4\tilde M(z,y_1)^2\hat{P}(z,y_2+y_3,y_1,y_2+y_3)\hat{P}(z,0,y_1,0)^3\right).
\end{align*}

\begin{center}
 	\begin{figure}[h!]
 		\begin{minipage}{1\textwidth}
 			\begin{center}
 				{\includegraphics[width=0.9\textwidth]{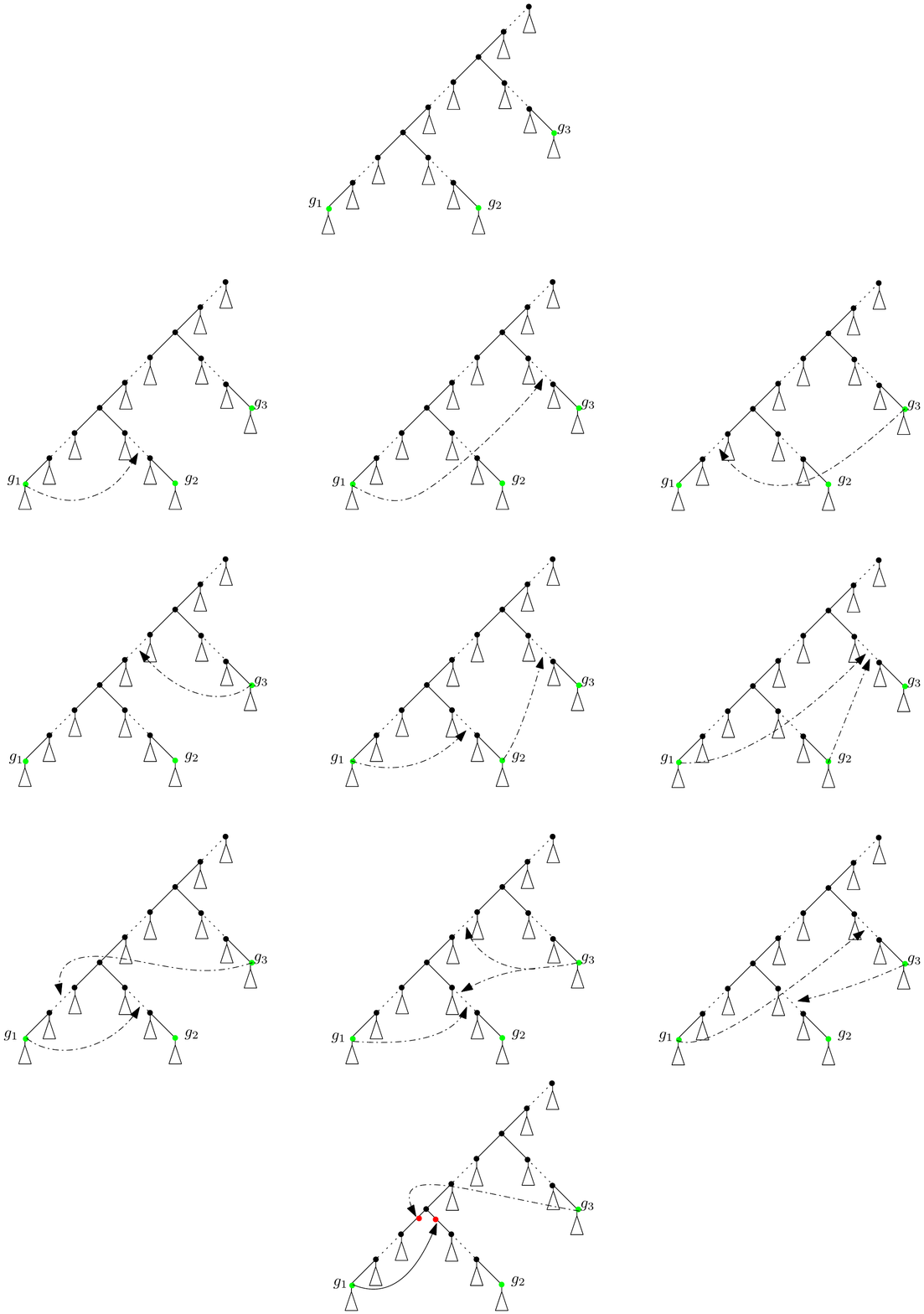}}
 			\end{center}
 		\end{minipage}
 		\caption{The colored Motzkin skeletons arising from the sparsened skeleton in
Figure~\ref{spsk-k3}-(d) classified according to the indicated pointing rules; the missing pointers from green vertices have no restrictions except that they are not allowed to point at vertices on a path. The exponential generating function of all except the last have to be added, whereas the one from the last one has to be subtracted.}
 		\label{T3-4}
 	\end{figure}
 \end{center}

Finally, the cases for the colored Motzkin skeletons arising from the sparsened skeleton in
Figure~\ref{spsk-k3}-(d) are classified according to the indicated pointing rules in Figure~\ref{T3-4}. Here, the exponential generating function of all except the last one have to be added, whereas the exponential generating function of the last one has to be subtracted. This yields
\begin{align*}
\rcp{2}&\textbf{Y}\left(z^5\tilde{M}(z,y_1+y_2+y_3)^3\hat{P}(z,0,y_1+y_2+y_3,0)^5\right)\\
&+\textbf{Y}\left(z^5\tilde{M}(z,y_2+y_3)^3\hat{P}(z,y_1,y_2+y_3, y_1)\hat{P}(z,0,y_2+y_3,0)^4\right)\\
&+\textbf{Y}\left(z^5\tilde{M}(z,y_2+y_3)^3\hat{P}(z,y_1, y_2+y_3,y_1)\hat{P}(z,0,y_2+y_3,0)^4\right)\\
&+\textbf{Y}\left(z^5\tilde{M}(z,y_1+y_2)^3\hat{P}(z,y_3,y_1+y_2,y_3)\hat{P}(z,0,y_1+y_2,0)^4\right)\\
&+\rcp{2}\textbf{Y}\left(z^5 \tilde{M}(z,y_1+y_2)^3\hat{P}(z,y_3,y_1+y_2,y_3)\hat{P}(z,0,y_1+y_2,0)^4\right)\\
&+\textbf{Y}\left(z^5\tilde{M}(z,y_3)^3\hat{P}(z,y_1,y_3,y_1)\hat{P}(z,y_2,y_3,y_2)\hat{P}(z,0,y_3,0)^3\right)\\
&+\rcp{2}\textbf{Y}\left(z^5\tilde{M}(z,y_3)^3\hat{P}(z,y_1+y_2,y_3,y_1+y_2)\hat{P}(z,0,y_3,0)^4\right)\\
&+\textbf{Y}\left(z^5\tilde{M}(z,y_2)^3\hat{P}(z,y_1,y_2,y_1)\hat{P}(z,y_3,y_2,y_3)\hat{P}(z,0,y_2,0)^3\right)\\
&+\textbf{Y}\left(z^5\tilde{M}(z,y_2)^3\hat{P}(z,y_1+y_3,y_2,y_1+y_3)\hat{P}(z,y_3,y_2,y_3)\hat{P}(z,0,y_2,0)^3\right)\\
&+\textbf{Y}\left(z^5\tilde{M}(z,y_2)^3\hat{P}(z,y_1,y_2,y_1)\hat{P}(z,y_3,y_2,y_3)\hat{P}(z,0,y_2,0)^3\right)\\
&-\textbf{Y}\left(z^5\tilde{M}(z,y_2)^3\hat{P}(z,0,y_2,y_1)\hat{P}(z,0,y_2,y_3)\hat{P}(z,0,y_2,0)^3\right).
\end{align*}

Adding up the above three exponential generating functions corresponding to the sparsened
skeletons in Figure~\ref{spsk-k3}-(b), (c), (d), then adding to this sum the one arising from the
sparsened skeleton in Figure~\ref{spsk-k3}-(a), and finally dividing by $2^k=8$ (since every
tree-child networks is generated from this $8$ times), we get
\[
T_{3}(z)=z\frac{-35z^{12}+175z^{10}-(34z^{12}+175z^{10})\sqrt{1-2z^2}}{(1-2z^2)^{11/2}}.
\]
Consequently,
\begin{align*}
T_{3,2n+1}=(2n+1)!\Bigg(&\frac{n^2(n-1)(n-2)(n-3)(n-4)}{3(2n-1)2^{n}}\binom{2n}{n}\\
&\qquad-\frac{2^{n}}{192}(n-2)(n-3)(n-4)(48n-65)\Bigg)
\end{align*}
and
\begin{align}
\tilde{T}_{3,\ell}=\ell!\Bigg(&\frac{(\ell+2)^2(\ell+1)\ell(\ell-1)(\ell-2)}{12(2\ell+3)2^{\ell}}\binom{2\ell+4}{\ell+2}
-\frac{2^{\ell}}{48}\ell(\ell-1)(\ell-2)(48\ell+31)\Bigg),
\label{eq:TL3}
\end{align}
The latter sequence (for $\ell\geq 4$) starts with
\[
2544, 154500,  6494400, 241204950, 8609378400,\ldots
\]
which is in accordance with the values in Table 3 of \cite{Zh2}.

\section{Normal Networks with $2$ and $3$ Reticulation Vertices}\label{normal}

Here, we explain how to modify our approach from \cite{FuGiMa} to get the correct exponential generating functions for the number of normal networks with $k=2$ and $k=3$. As explained in the last paragraph in Section~\ref{meth-FuGiMa}, apart from avoiding to create cycles in the pointing step, we also have to be careful not to create the two near-cycles discussed from that paragraph.

In fact, avoiding the creation of cycles and the near-cycles in Figure~\ref{cycle}-(b) is done by considering the same cases as in the last section. Then, we will subtract all the networks which contain the near-cycles in Figure~\ref{cycle}-(c). For technical reasons, it will be advantageous to replace the exponential generating function (\ref{path-nn}) for paths for normal networks  by the following more detailed one:
\[
\tilde{P}(z,y,\tilde{y},\bar{y},\hat{y})=\frac{1+z\hat y}{1-(z+z^2y+z^2\bar{y})\tilde M(z,\tilde y)},
\]
where the only difference to the previous one is that now $y$ marks unary vertices on the path and $\bar{y}$ marks unary vertices which are children of vertices on the path (in (\ref{path-nn}) both of these vertices were marked by $y$).

\medskip
\subsection{Normal Networks with Two Reticulation Vertices} We again start from the two sparsened skeletons in Figure~\ref{spsk-k2}.

\begin{center}
	\begin{figure}
		\begin{minipage}{1\textwidth}
			\begin{center}
				{\includegraphics[width=0.8\textwidth]{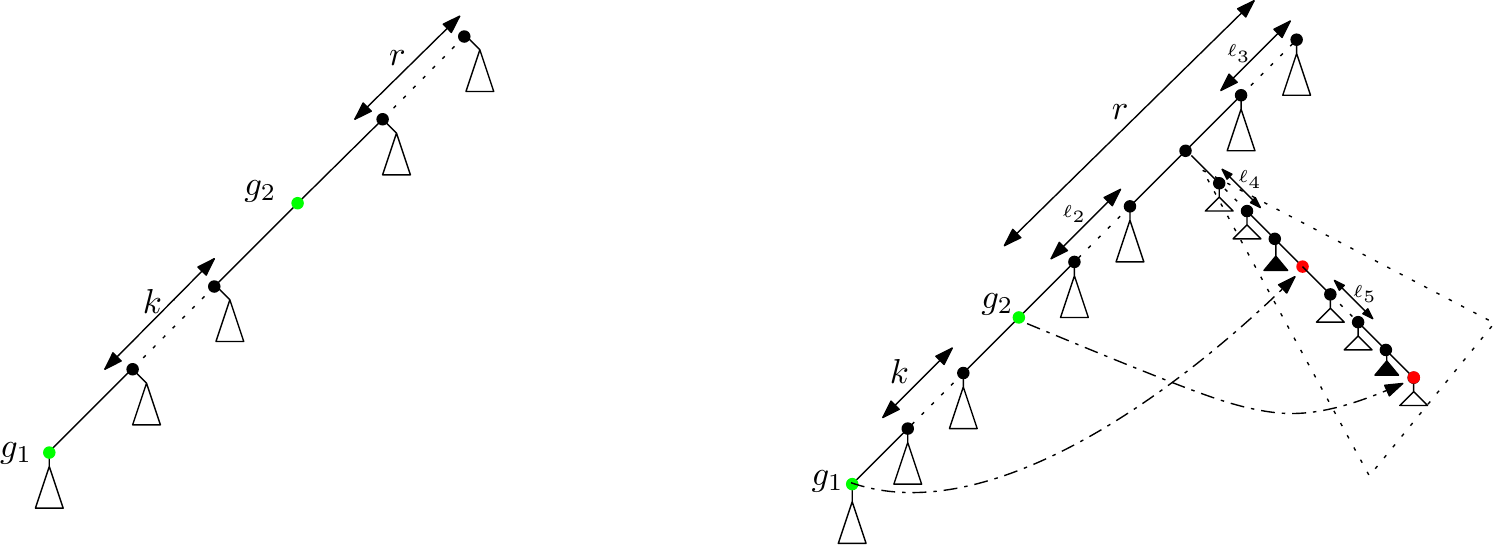}}
			\end{center}
		\end{minipage}
		\caption{The colored Motzkin skeletons arising from Figure~\ref{spsk-k2}-(a) and
the networks containing near-cycles which have to be subtracted. Solid subtrees mean that these
subtrees must be there, because otherwise the tree-child condition would be violated.}
		\label{N2-1}
	\end{figure}
\end{center}

This time, we also have to consider the sparsened skeleton in Figure~\ref{spsk-k2}-(a), because it
is possible to create the near-cycles in Figure~\ref{cycle}-(c) (which have to be subtracted). We
consider in Figure~\ref{N2-1} the colored Motzkin skeletons which arise from that sparsened
skeleton (left) and the networks containing the near-cycles in Figure~\ref{cycle}-(c) which have
to be subtracted (right; solid subtrees mean that these subtrees must be there; also the pointing rules are indicated). Thus the
exponential generating function must satisfy the equation
\[
\partial_{y_1}\partial_{y_2}\left(z^2\tilde{M}(z,0)\tilde{P}(z,0,y_1,0,0)\tilde{P}(z,0,y_1+y_2,0,0)\right)\Big\vert_{y_1=y_2=0}-z^7\tilde{M}(z,0)^4\tilde{P}(z,0,0,0,0)^5,
\]
where in the first term, the $z^2$ counts the two green vertices, $\tilde{M}(z,0)$ counts the tree dangling from $g_1$, and $\tilde{P}(z,0,y_1,0,0)$ and $\tilde{P}(z,0,y_1+y_2,0,0)$ count the two paths $k$ and $r$. In the second term, which is the subtraction term, $z^7$ counts the two green vertices, the two endpoints of the pointers, the last common ancestor of these four vertices, and the roots of the two solid subtrees in the network on the right of Figure~\ref{N2-1}; $\tilde{M}(z,0)^4$ counts the subtree dangling from $g_1$, the subtree before the vertex to which the pointer from $g_1$ points, and the two subtrees before and after the vertex to which the pointer of $g_2$ points; and $\tilde{P}(z,0,0,0,0)^5$ counts the five paths $k, \ell_2, \ell_3, \ell_4, \ell_5$.

For the sparsened skeleton in Figure~\ref{spsk-k2}-(b), we use the colored Motzkin skeletons from
Figure~\ref{T2-2}, where in the skeletons from Figure~\ref{T2-2}-(ii), the pointer of $g_2$ is
neither allowed to point at a vertex on a path nor at the child of a vertex on a path. This
implies that the exponential generating function satisfies the equation
\begin{align*}
\frac{1}{2}\partial_{y_1}\partial_{y_2}&\left(z^3\tilde{M}(z,y_1)\tilde{M}(z,y_2)\tilde{P}(z,0,y_1+y_2,y_1,0)\tilde{P}(z,0,y_1+y_2,y_2,0)
\tilde{P}(z,0,y_1+y_2,0,0)\right.\\
&\left.\qquad+z^3\tilde{M}(z,y_2)\tilde{M}(z,0)\tilde{P}(z,y_1,y_2,0,0)\tilde{P}(z,0,y_2,0,0)^2\right)\Big\vert_{y_1=y_2=0}.
\end{align*}

Now, adding the exponential generating functions of the above two cases and dividing by $2^k=4$ (since every network is obtained from this procedure exactly $4$ times) gives
\[
N_{2}(z)=z\frac{11z^8-66z^6+50z^4-8z^2-(-28z^6+42z^4-8z^2)\sqrt{1-2z^2}}{(1-2z^2)^{7/2}}
\]
from which we have
\[
N_{2,2n+1}=(2n+1)!\left(\frac{n(3n-7)(n^2+9n-4)}{(2n-1)2^n}\binom{2n}{n}-2^{n-1}(n+1)(3n-7)\right)
\]
and
\begin{equation}
\tilde{N}_{2,\ell}=\ell!\left(\frac{(\ell+1)(3\ell-4)(\ell^2+11\ell+6)}{6(2\ell+1)2^{\ell}}\binom{2\ell+2}{\ell+1}-2^{\ell}(\ell+2)(3\ell-4)\right).
\label{eq:NL2}
\end{equation}
The latter sequence starts with (for $\ell\geq 4$)
\[
48, 2310, 78120,  2377620, 70749000,\ldots
\]
which is in accordance with the values from Table 1 in \cite{Zh1}.

\medskip
\subsection{Normal Networks with Three Reticulation Vertices} Finally, we consider normal networks with three reticulation vertices. The sparsened skeletons are again in Figure~\ref{spsk-k3}. We will consider below the exponential generating function for the colored Motzkin skeletons arising from each of them.

\begin{center}
	\begin{figure}
		\begin{minipage}{1\textwidth}
			\begin{center}
				{\includegraphics[width=1\textwidth]{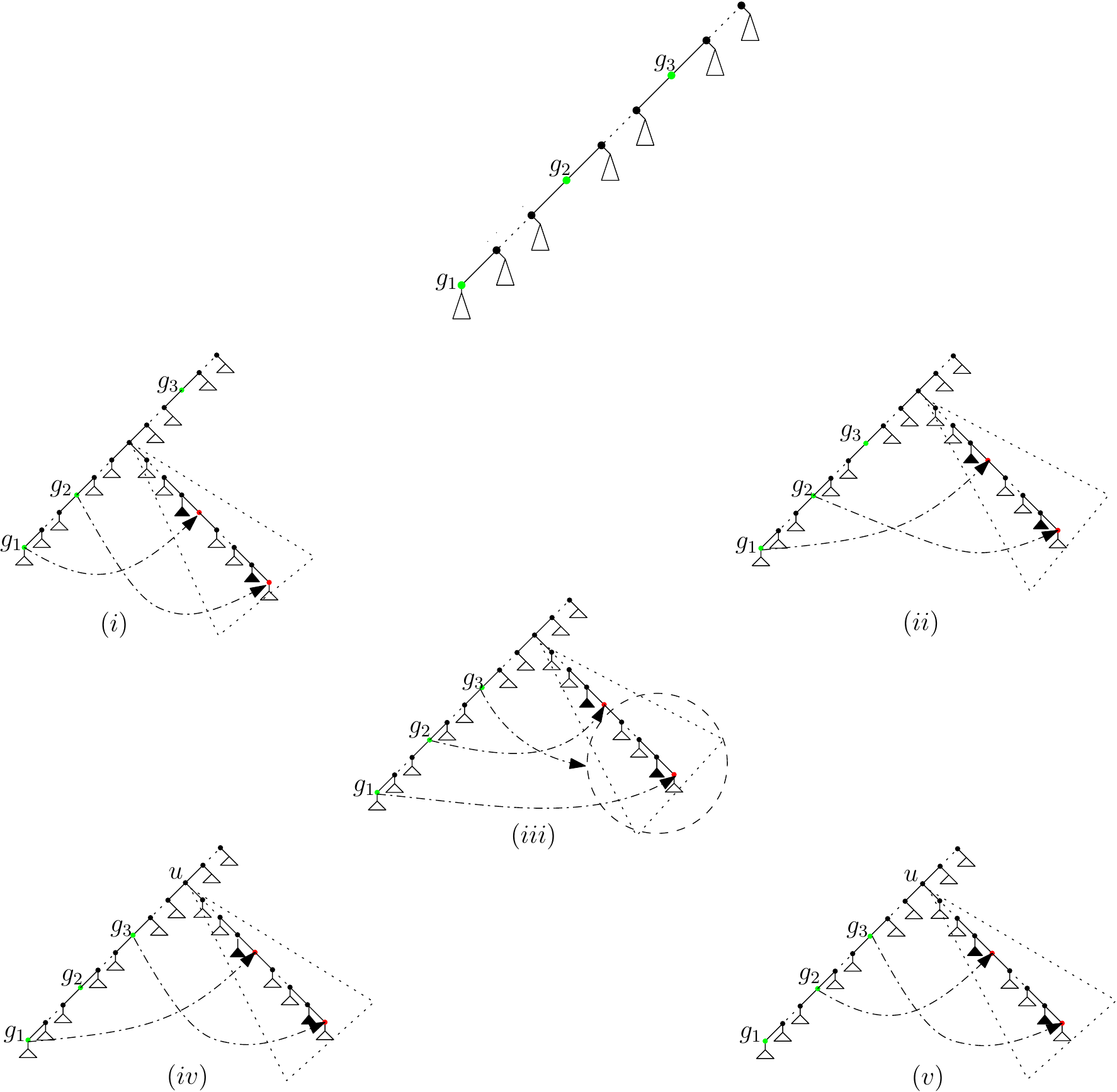}}
			\end{center}
		\end{minipage}
		\caption{The colored Motzkin skeletons arising from the sparsened skeleton in
Figure~\ref{spsk-k3}-(a) with the five subtraction cases due to the creation of near-cycles from
Figure~\ref{cycle}-(c). The pointing rules for $g_3$ in the five cases are as follows: (i): $g_3$ can point at any vertex so that the normal condition is not violated; (ii) and (iii): these are the cases where $g_1$ and $g_2$ point ar a vertex on the same path. Thus,
in (iii), $g_3$ must point to a vertex inside the circle, because otherwise there is no near-cycle
from Figure~\ref{cycle}-(c). (iv) and (v): the remaining two cases with $g_1$ and $g_2$ not
pointing at vertices on the same path. Thus, in (iv) resp. (v), $g_2$ resp. $g_1$ is not allowed
to point at a vertex of path between $u$ and the endpoint of the pointer from $g_1$ resp. $g_2$ and no vertex after that endpoint.}
		\label{N3-1}
	\end{figure}
\end{center}

First, for the sparsened skeleton in Figure~\ref{spsk-k3}-(a), we consider the colored
Motzkin-skeletons from Figure~\ref{N3-1} with the networks which have to be subtracted because they contain the near-cycle from Figure~\ref{cycle}-(c) (if the pointing of the missing pointers in the subtraction cases is not merely restricted by the normal condition, we explain it in the figure; solid subtrees mean again that they must be there). Overall, we get

\begin{align*}
\textbf{Y}&\left(z^3M(z,0)\tilde{P}(z,0,y_1,0,0)\tilde{P}(z,0,y_1+y_2,0,0)\tilde{P}(z,0,y_1+y_2+y_3,0,0)\right)\\
&-\partial y_3\left(z^8\tilde{M}(z,0)^4\tilde{P}(z,0,0,0,0)^5\tilde{P}(z,0,y_3,0,0)\right)\Big\vert_{y_1=0}\\
&-\partial y_3\left(z^8\tilde{M}(z,0)\tilde{M}(z,y_3)^3\tilde{P}(z,y_3, y_3,y_3,0)^2\tilde{P}(z,0,0,0,0)^2\tilde{P}(z,0,y_3,0,0)^2\right)\Big\vert_{y_1=0}\\
&-\partial y_3\left(z^8\tilde{M}(z,0)^2\tilde{M}(z,y_3)^2\tilde{P}(z,y_3,y_3,y_3,0)\tilde{P}(z,0,0,0,0)^5\right)\Big\vert_{y_1=0}\\
&-\partial y_2\left(z^8\tilde{M}(z,0)^3\tilde{M}(z,y_2)\tilde{P}(z,0,y_2,y_2,0)\tilde{P}(z,0,0,0,0)^2\tilde{P}(z,0,y_2,0,0)^3\right)\Big\vert_{y_1=0}\\
&-\partial y_1\left(z^8\tilde{M}(z,0)^3\tilde{M}(z,y_1)\tilde{P}(z,0,y_1,y_1,0)\tilde{P}(z,0,0,0,0)\tilde{P}(z,0,y_1,0,0)^4\right)\Big\vert_{y_1=0},
\end{align*}
where ${\bf Y}(\cdot)$ is as in the last section and the last five terms correspond to the cases (i) until (iv) in Figure~\ref{N3-1} in that order.

\begin{center}
	\begin{figure}
		\begin{minipage}{1\textwidth}
			\begin{center}
				{\includegraphics[width=.8\textwidth]{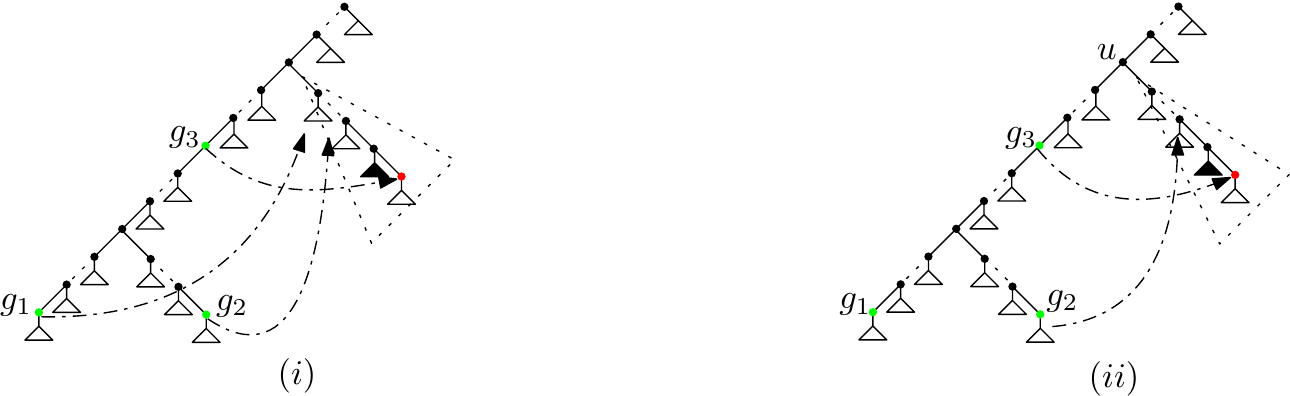}}
			\end{center}
		\end{minipage}
		\caption{The subtraction cases for the colored Motzkin skeletons arising from the
sparsened skeleton in Figure~\ref{spsk-k3}-(b). In (i), the pointers of $g_1,g_2$ both cross the
pointer of $g_3$ and point at vertices on the same path, whereas in (ii) only the pointers of
$g_2$ and $g_3$ cross and point at vertices on the same path, \emph{i.e.}, $g_1$ is not allowed to point at a vertex of the path between $u$ and the endpoint of the pointer from $g_3$.}
		\label{N3-2}
	\end{figure}
\end{center}

Next, for the sparsened skeleton in Figure~\ref{spsk-k3}-(b), we consider first the cases from
Figure~\ref{T3-2} to create networks which respect the tree-child condition and do not contain the
near-cycle in Figure~\ref{cycle}-(b) (for this, we do not need the third network in
Figure~\ref{T3-2}, because its creation will be avoided by our method). Then, we subtract all networks containing the near-cycle in Figure~\ref{cycle}-(c); see Figure~\ref{N3-2} where all the networks we have to subtract are listed (and restrictions to the pointing of the missing pointers is explained in case pointing is not merely restricted by the normal condition; solid subtrees meant that they must be there). This gives
\begin{align*}
\rcp{2}&\textbf{Y}\left(z^4\tilde{M}(z,y_1)\tilde{M}(z,y_2)\tilde{P}(z,0,y_1+y_2,y_2,0)P(z,0,y_1+y_2,y_1,0)\tilde{P}(z,0,y_1+y_2,0,0)\right.\\
&\qquad\left.\times\tilde{P}(z,0,y_1+y_2+y_3,0,0)\right)\\
&+\textbf{Y}\left(z^4\tilde{M}(z,0)\tilde{M}(z,y_2)\tilde{P}(z,y_1,y_2,0,0)\tilde{P}(z,0,y_2,0,0)^2
\tilde{P}(z,0,y_2+y_3,0,0)\right)\\
&-\rcp{2}\partial y_1 \partial y_2\left(z^7\tilde{M}(z,0)^4\tilde{P}(z,y_1+y_2,0,0,0)\tilde{P}(z,0,0,0,0)^5\right)\Big\vert_{y_1=y_2=0}\\
&-\partial y_1\partial y_2\left(z^7\tilde{M}(z,0)\tilde{M}(z,y_1)^3\tilde{P}(z,y_2,y_1,y_1,0)\tilde{P}(z,y_1,y_1,y_1,0)\tilde{P}(z,0,y_1,0,0)^4\right)
\Big\vert_{y_1=y_2=0},
\end{align*}
where the last two terms correspond to the cases (i) and (ii) in Figure~\ref{N3-2} in that order.

\begin{center}
	\begin{figure}
		\begin{minipage}{.8\textwidth}
			\begin{center}
				{\includegraphics[width=1.05\textwidth]{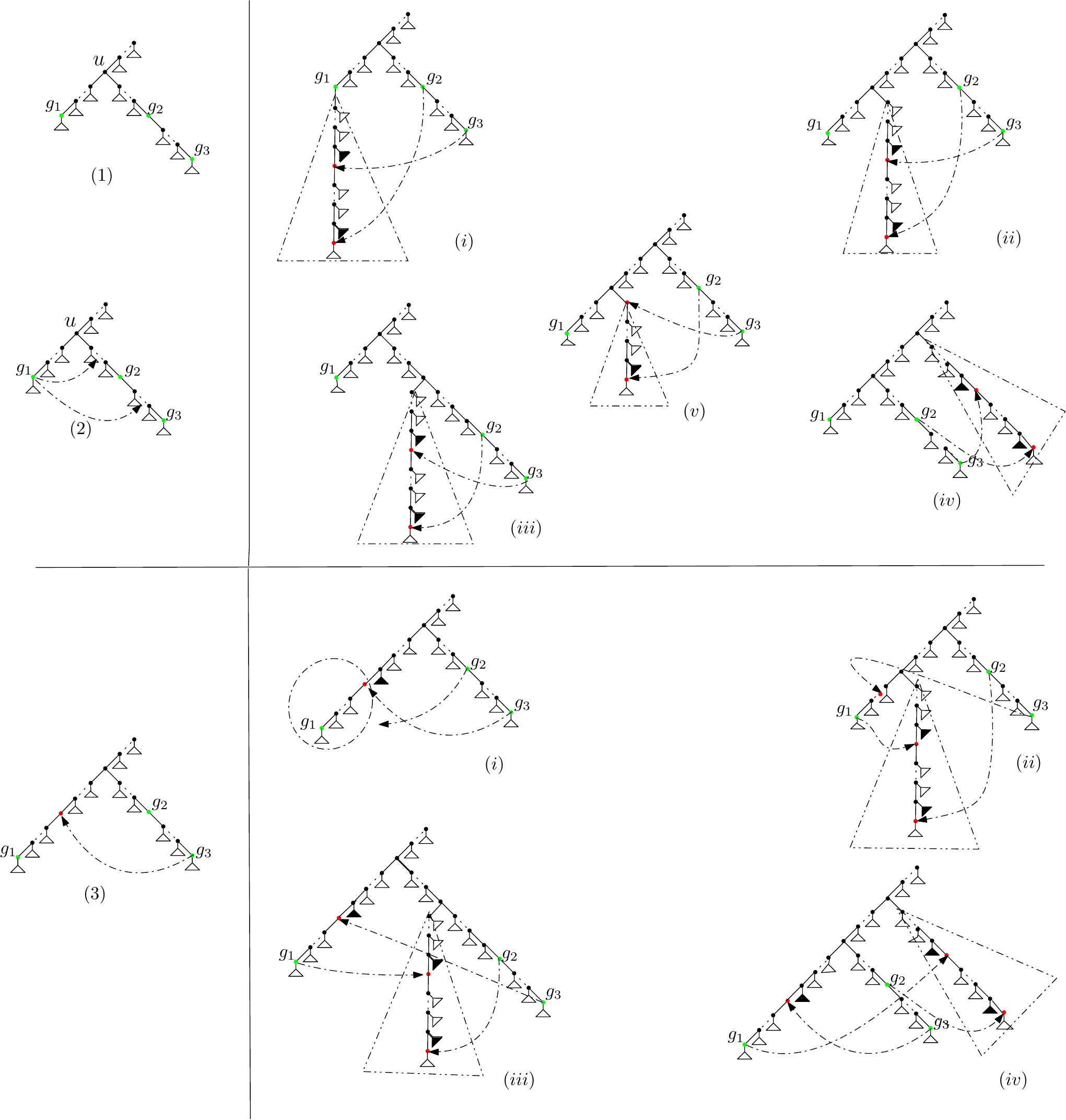}}
			\end{center}
		\end{minipage}
		\caption{The subtraction cases for the colored Motzkin skeletons arising from the
sparsened skeleton in Figure~\ref{spsk-k3}-(c). The first column contains the cases from
Figure~\ref{T3-3} which contain the near-cycles in the second column. In the first row, only $g_2$
and $g_3$ are on the same path, so a near-cycle is only created if their pointers cross and point
at a vertex on the same path. That path can be after $u$ (cases (i), (ii), and (v) in the second
column); between $u$ and $g_2$ (case (iii) in the second column); or before $u$ (case (iv) in the
second column). In the second row $g_2, g_3,$ and $g_1$ are all on the same path (due to the
pointer of $g_3$). So, one needs to subtract the cases where the pointers from $g_2$ and $g_3$
cross and point at vertices on the same path (case (i) in the second column) and where the
pointers from $g_1$ and $g_2$ cross and point at vertices on the same path (cases (i), (ii), (iii), and (iv) in the second column).}
		\label{N3-3}
	\end{figure}
\end{center}

Third, for the sparsened skeleton in Figure~\ref{spsk-k3}-(c), we start again with the cases from Figure~\ref{T3-3} and use them to create networks which do not contain the cycles from Figure~\ref{cycle}-(a) and the near-cycles from Figure~\ref{cycle}-(b). This gives
\begin{align*}
\textbf{Y}&\left(z^4\tilde{M}(z,y_1)\tilde{M}(z,y_2+y_3)\tilde{P}(z,0,y_1+y_2+y_3,y_2+y_3,0)\tilde{P}(z,0,y_1+y_3,y_1,0)\right.\\
&\qquad\left.\times\tilde{P}(z,0,y_1+y_2+y_3,y_1,0)\tilde{P}(0,y_1+y_2+y_3,0,0)\right)\\
&+\textbf{Y}\left(z^4\tilde{M}(z,0)\tilde{M}(z,y_2+y_3)\tilde{P}(z,y_1,y_2+y_3,0,0)\tilde{P}(z,0,y_2+y_3,0,0)^2\tilde{P}(z,0,y_3,0,0)\right)\\
&+\textbf{Y}\left(z^4\tilde{M}(z,0)\tilde{M}(z,y_2+y_3)\tilde{P}(z,y_1,y_3,0,0)\tilde{P}(z,0,y_2+y_3,y_2,0)\tilde{P}(z,0,y_2+y_3,0,0)^2\right)\\
&+\textbf{Y}\left(z^4\tilde{M}(z,y_1)\tilde{M}(z,y_2)\tilde{P}(z,y_3,y_1+y_2,y_2,0)\tilde{P}(z,0,y_1,0,0)\tilde{P}(z,0,y_1+y_2,0,0)^2\right)\\
&+\textbf{Y}\left(z^4\tilde{M}(z,y_1)\tilde{M}(z,y_3)\tilde{P}(z,y_2,y_1+y_3,y_3,0)\tilde{P}(z,0,y_1+y_3,y_1,0)\tilde{P}(z,0,y_1+y_3,0,0)^2\right)\\
&+\textbf{Y}\left(z^4\tilde{M}(z,0)\tilde{M}(z,y_1)\tilde{P}(z,y_2,y_3,0,0)\tilde{P}(z,y_1,y_3,0,0)\tilde{P}(z,0,y_3,0,0)^2\right)\\
&+\textbf{Y}\left(z^4\tilde{M}(z,0)\tilde{M}(z,y_1)\tilde{P}(z,y_2+y_3,y_1,0,0)\tilde{P}(z,0,y_1,0,0)^3\right).
\end{align*}
Then, we subtract from this the exponential generating functions of all cases where networks
contain the cycles from Figure~\ref{cycle}-(c). Here, in contrast to the other cases, there are
many such situations and all of them are listed in Figure~\ref{N3-3}. The sum of the exponential
generating functions of all these cases is given by

\begin{align*}
\partial {y_1}&\left(z^8\tilde{M}(z,0)^3\tilde{M}(z,y_1)\tilde{P}(z,y_1, y_1,y_1,0)^2\tilde{P}(z,0,y_1,0,0)^2\tilde{P}(z,0,0,0)^2\right.\\
&\qquad+\left.z^9\tilde{M}(z,0) \tilde{M}(z,y_1)^4\tilde{P}(z,y_1,y_1,y_1,0)^4\tilde{P}(z,0,y_1,0,0)^3\right.\\
&\qquad+\left.z^9\tilde{M}(z,0)\tilde{M}(z,y_1)^4\tilde{P}(z,y_1, y_1,y_1,0)^3\tilde{P}(z,y_1,y_1,y_1,y_1)^2\tilde{P}(z,0,y_1,0,0)^2\right.\\
&\qquad+\left.z^9\tilde{M}(z,0)\tilde{M}(z,y_1)^4\tilde{P}(z,y_1,y_1,y_1,0)^4\tilde{P}(z,0,y_1,0,0)^3\right.\\
&\qquad+\left.z^8\tilde{M}(z,0)\tilde{M}(z,y_1)^3\tilde{P}(z,0,y_1,y_1,0)^2\tilde{P}(z,y_1, y_1, y_1,0)\tilde{P}(z,0,y_1,0,0)^3\right)\Big\vert_{y_1=0}\\
&+\partial{y_1}\partial{y_2}\left(z^6\tilde{M}(z,y_2)\tilde{M}(z,y_1)^2\tilde{P}(z,y_2,y_1+y_2,y_2,0)\tilde{P}(z,0,y_1,0,0)^4\right)\Big\vert_{y_1=0, y_2=0}\\
&+z^{10}\tilde{M}(z,0)^5\tilde{P}(z,0,0,0)^8+z^{11}\tilde{M}(z,0)^6\tilde{P}(z,0,0,0,0)^8+z^{11}\tilde{M}(z,0)^6\tilde{P}(z,0,0,0,0)^8,
\end{align*}
where the terms in the first bracket correspond to the networks from row one and column two of Figure~\ref{N3-3} in the order from (i) to (v) and the remaining terms correspond to the networks from row two and column two of Figure~\ref{N3-3} in the order from (i) to (iv).

Finally, we consider the colored Motzkin skeletons arising from the sparsened skeleton in Figure~\ref{spsk-k3}-(d). Here, the creation of the near-cycles from Figure~\ref{cycle}-(c) is impossible. Thus, we only have to consider the cases from Figure~\ref{T3-4} (except the last case in this figure since the occurrence of these networks is already ruled out by the normal condition) and make sure that the creation of the cycles and near-cycles from Figure~\ref{cycle}-(a) and Figure~\ref{cycle}-(b) is avoided. Overall, we obtain

\begin{align*}
\rcp{2}\textbf{Y}&\left(z^5\tilde{M}(z,y_1+y_2)\tilde{M}(z,y_1+y_3)\tilde{M}(z,y_2+y_3)\tilde{P}(z,0,y_1+y_2+y_3,y_1+y_2,0)\right.\\
&\qquad\times\left.\tilde{P}(z,0,y_1+y_2+y_3,y_1+y_3,0)\tilde{P}(z,0,y_1+y_2+y_3, y_2+y_3,0)\right.\\
&\qquad\times\left.\tilde{P}(z,0,y_1+y_2+y_3,y_3,0)\tilde{P}(z,0,y_1+y_2+y_3,0,0)\right)\\
&+\textbf{Y}\left(z^5\tilde{M}(z,y_2)\tilde{M}(z,y_3)\tilde{M}(z,y_2+y_3)\tilde{P}(z,y_1,y_2+y_3, y_3,0)\right.\\
&\qquad\times\left.\tilde{P}(z,0,y_2+y_3,y_2,0)\tilde{P}(z,0,y_2+y_3,y_3,0)^2\tilde{P}(z,0,y_2+y_3,0,0)\right)\\
&+\textbf{Y}\left(z^5\tilde{M}(z,y_2)\tilde{M}(z,y_3)\tilde{M}(z,y_2+y_3)\tilde{P}(z,y_1,y_2+y_3,y_2,0)\tilde{P}(z,0,y_2+y_3,y_3,0)\right.\\
&\qquad\times\left.\tilde{P}(z,0,y_2+y_3,y_2,0)\tilde{P}(z,0,y_2+y_3,0,0)^2\right)\\
&+\textbf{Y}\left(z^5\tilde{M}(z,y_1)\tilde{M}(z,y_2)\tilde{M}(z,y_1+y_2)\tilde{P}(z,y_3,y_1+y_2,y_2,y_3)\tilde{P}(z,0,y_1+y_2,y_1,0)\right.\\
&\qquad\times\left.\tilde{P}(z,0,y_1+y_2,y_2,0)\tilde{P}(z,0,y_1+y_2,0,0)^2\right)\\
&+\frac{1}{2}\textbf{Y}\left(z^5\tilde{M}(z,y_1)\tilde{M}(z,y_2)\tilde{M}(z,y_1+y_2)\tilde{P}(z,y_3,y_1+y_2,0,0)\tilde{P}(z,0,y_1+y_2,y_1,0)\right.\\
&\qquad\times\left.\tilde{P}(z,0,y_1+y_2,y_2,0)\tilde{P}(z,0,y_1+y_2,0,0)^2\right)\\
&+\textbf{Y}\left(z^5\tilde{M}(z,y_3)^2\tilde{M}(z,0)\tilde{P}(z,y_1,y_3,0,0)\tilde{P}(z,y_2,y_3,0,0)\tilde{P}(z,0,y_3,0,0)^3\right)\\
&+\rcp{2}\textbf{Y}\left(z^5\tilde{M}(z,y_3)^2\tilde{M}(z,0)\tilde{P}(z,y_1+y_2,y_3,0,0)\tilde{P}(z,0,y_3,0,0)^4\right)\\
&+\textbf{Y}\left(z^5\tilde{M}(z,y_2)^2\tilde{M}(z,0)\tilde{P}(z,y_1,y_2,0,0)\tilde{P}(z,y_3,y_2,0,y_3)\tilde{P}(z,0,y_2,0,0)^3\right)\\
&+\textbf{Y}\left(z^5\tilde{M}(z,y_2)^2\tilde{M}(z,0)\tilde{P}(z,y_1+y_3,y_2,0,y_3)\tilde{P}(z,y_3,y_2,0,0)\tilde{P}(z,0,y_2,0,0)^3\right)\\
&+\textbf{Y}\left(z^5\tilde{M}(z,y_2)^2\tilde{M}(z,0)\tilde{P}(z,y_1,y_2,0,0)\tilde{P}(z,y_3,y_2,0,0)\tilde{P}(z,0,y_2,0,0)^3\right).
\end{align*}

Now, by combining all the contributions above and dividing the result by $2^k=8$ (since every normal network is created by the above procedure exactly $8$ times), we have
\begin{align*}
N_3(z)=z&\left(\frac{877z^{12}-3065z^{10}+2392z^8-628z^6+64z^4}{(1-2z^2)^{11/2}}\right.\\
&\qquad\left.-\frac{110z^{12}-1455z^{10}+1860z^{8}-564z^6+64z^4}{(1-2z^2)^5}\right).
\end{align*}
From this, we obtain that
\begin{align*}
N_{3,2n+1}=(2n+1)!\Bigg(&\frac{n(n-1)(n^4+15n^3-158n^2+324n+40)}{3(2n-1)2^{n}}\binom{2n}{n}\\
&\qquad-\frac{2^{n}}{192}(144n^4-751n^3-1089n^2+9106n-7080)\Bigg)
\end{align*}
and
\begin{align}
\tilde{N}_{3,\ell}=\ell!\Bigg(&\frac{(\ell+2)(\ell+1)(\ell^4+23\ell^3-44\ell^2-96\ell+192)}{12(2\ell+3)2^{\ell}}\binom{2\ell+4}{\ell+2} \nonumber  \\
&\qquad-\frac{2^{\ell}}{48}(144\ell^4+401\ell^3-2139\ell^2+346\ell+3072)\Bigg),
\label{eq:NL3}
\end{align}
The latter sequence (for $\ell\geq 5$) starts with
\[
1920, 184680, 11059650, 547444800,\ldots.
\]
The first two values coincide with the values given in Table~1 in \cite{Zh1}. However, the next
two are different from the ones erroneously given in \cite{Zh1} as $11038530$ and $536524830$. In
private communication, L. Zhang told us that these wrong values in \cite{Zh1} resulted from an
overflow problem in his C++ program and that our values are indeed the correct ones. In fact, it
can be verified that the previous values are erroneous, because
\[
\frac{11038530}{7!}=\frac{367951}{168}=\frac{367951}{2^3\cdot 3\cdot 7}
\]
and
\[
\frac{536524830}{8!}=\frac{5961387}{448}=\frac{5961387}{2^6\cdot 7}
\]
which is impossible since the denominators have to be powers of $2$; see Corollary~\ref{cor-app} in the Appendix.

\section*{Acknowledgment}

We thank the two anonymous reviewers for a careful reading of the manuscript and many constructive remarks.

\section*{Appendix}

In this appendix, we want to find the answer to the following question:

\medskip
\noindent{\bf Question:}\ {\it Given a phylogenetic network $N$, how many different leaf-labeled networks can be generated from $N$ by labeling its leaves?}

\begin{figure}[h!]
\begin{center}
\begin{tikzpicture}[scale=0.7]
\draw (0cm,0cm) node[circle,inner sep=0.8,fill,draw] (1) {};
\draw (-0.8cm,-0.8cm) node[circle,inner sep=0.8,fill,draw] (2) {};
\draw (0.8cm,-0.8cm) node[circle,inner sep=0.8,fill,draw] (3) {};
\draw (0cm,-1.6cm) node[circle,inner sep=0.8,fill,draw] (4) {};
\draw (0cm,-2.4cm) node[circle,inner sep=0.8,fill,draw] (5) {};
\draw (0cm,-3.2cm) node[circle,inner sep=0.8,fill,draw] (6) {};
\draw (0cm,-4cm) node[circle,inner sep=0.8,fill,color=green!60,draw] (7) {};
\draw (-0.5cm,-4.8cm) node[circle,inner sep=0.8,fill,draw] (8) {};
\draw (0.5cm,-4.8cm) node[circle,inner sep=0.8,fill,draw] (9) {};
\draw (-1.4cm,-2.4cm) node (10) {(a)};
\draw (0cm,0.5cm) node (11) {};
\draw (0cm,-5.3cm) node (12) {};

\draw (1)--(2); \draw (1)--(3); \draw (2)--(4); \draw (3)--(4);
\draw (4)--(5); \draw (6)--(7); \draw (7)--(8); \draw (7)--(9);
\draw (3) to [out=-90,in=45] (6); \draw (2) to [out=-90,in=135] (6);
\draw[dashed] (11)--(12);

\draw (5cm,0cm) node[circle,inner sep=0.8,fill,color=green!60,draw] (1) {};
\draw (4.2cm,-0.8cm) node[circle,inner sep=0.8,fill,draw] (2) {};
\draw (5.8cm,-0.8cm) node[circle,inner sep=0.8,fill,draw] (3) {};
\draw (5cm,-1.6cm) node[circle,inner sep=0.8,fill,draw] (4) {};
\draw (5cm,-2.4cm) node[circle,inner sep=0.8,fill,draw] (5) {};
\draw (5cm,-3.2cm) node[circle,inner sep=0.8,fill,draw] (6) {};
\draw (5cm,-4cm) node[circle,inner sep=0.8,fill,color=green!60,draw] (7) {};
\draw (4.5cm,-4.8cm) node[circle,inner sep=0.8,fill,draw] (8) {};
\draw (5.5cm,-4.8cm) node[circle,inner sep=0.8,fill,draw] (9) {};
\draw (4.55cm,-0.45cm) node[circle,inner sep=0.8,fill,draw] (10) {};
\draw (5.45cm,-0.45cm) node[circle,inner sep=0.8,fill,draw] (11) {};
\draw (4.8cm,-0.7cm) node[circle,inner sep=0.8,fill,draw] (12) {};
\draw (5.2cm,-0.7cm) node[circle,inner sep=0.8,fill,draw] (13) {};
\draw (3.6cm,-2.4cm) node (14) {(b)};
\draw (5cm,0.5cm) node (15) {};
\draw (5cm,-5.3cm) node (16) {};

\draw (1)--(2); \draw (1)--(3); \draw (2)--(4); \draw (3)--(4);
\draw (4)--(5); \draw (6)--(7); \draw (7)--(8); \draw (7)--(9);
\draw (3) to [out=-90,in=45] (6); \draw (2) to [out=-90,in=135] (6);
\draw (10)--(12); \draw (11)--(13);
\draw[dashed] (15)--(16);

\draw (10cm,-0.4cm) node[circle,inner sep=0.8,fill,color=green!60,draw] (1) {};
\draw (9.2cm,-1.2cm) node[circle,inner sep=0.8,fill,color=red!80,draw] (2) {};
\draw (10.8cm,-1.2cm) node[circle,inner sep=0.8,fill,color=red!80,draw] (3) {};
\draw (10cm,-2cm) node[circle,inner sep=0.8,fill,draw] (4) {};
\draw (10cm,-2.8cm) node[circle,inner sep=0.8,fill,draw] (5) {};
\draw (10cm,-3.6cm) node[circle,inner sep=0.8,fill,draw] (6) {};
\draw (10cm,-4.4cm) node[circle,inner sep=0.8,fill,draw] (7) {};
\draw (9.55cm,-0.85cm) node[circle,inner sep=0.8,fill,draw] (10) {};
\draw (10.45cm,-0.85cm) node[circle,inner sep=0.8,fill,draw] (11) {};
\draw (9.8cm,-1.1cm) node[circle,inner sep=0.8,fill,draw] (12) {};
\draw (10.2cm,-1.1cm) node[circle,inner sep=0.8,fill,draw] (13) {};
\draw (8.6cm,-2.4cm) node (14) {(c)};

\draw (1)--(2); \draw (1)--(3); \draw (2)--(4); \draw (3)--(4);
\draw (4)--(5); \draw (6)--(7);
\draw (3) to [out=-90,in=45] (6); \draw (2) to [out=-90,in=135] (6);
\draw (10)--(12); \draw (11)--(13);
\end{tikzpicture}
\end{center}
\caption{Three examples of phylogenetic networks with independent vertices in green and the only
independent pair of vertices in red. Note that in network (a), the root is symmetric but not
independent since all symmetric vertices are on the symmetry line (dashed line) as is the sole
leaf which is not contained in a proper subnetwork of a symmetric vertex. On the other hand,
in network (b), the root is independent, because there are two leaves not on the symmetry line (again drawn as a dashed line)
which are not in a proper subnetwork of symmetric vertices. Finally, network (c) is the only network
containing a symmetric pair since what is dangling from the two red vertices (namely, twice a reticulation vertex followed by a leaf)
is clearly identical. The values of $f$ for the three networks are $f=1$ (a), $f=2$ (b), and $f=2$ (c).}\label{fig-app}
\end{figure}
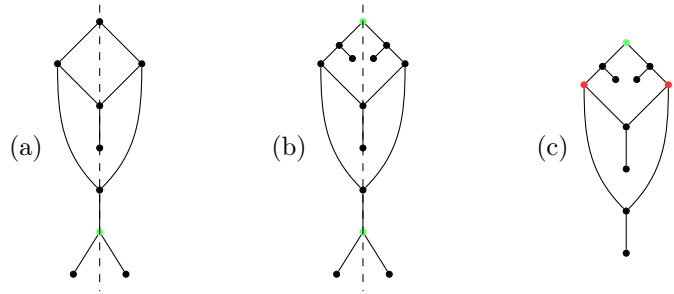

\medskip
For instance, in Example (a) in Figure~\ref{fig-app}, the answer to the question is $3$ since there are $3!=6$ possible labelings of the $3$ leaves, however, the leaf-labeled networks with the two lowest leaves are the same when the labels of the leaves are interchanged.

Note that for phylogenetic trees (which are leaf-labeled phylogenetic networks without reticulation vertices), the answer to the above question is known; see \cite{HeLiPe} and Section 2.4 in \cite{SeSt}.

Now for general phylogenetic networks, we denote the set of leaf-labeled networks from the above question by $P(N)$. Let $\tilde{N}\in P(N)$ be any network from $P(N)$, i.e., $N$ together with any labeling of its leaves. Then, by the Burnside lemma, we have
\[
\vert P(N)\vert=\frac{\ell!}{\vert F(\tilde{N})\vert},
\]
where $\ell$ is the number of labels of $\tilde{N}$ and $F(\tilde{N})$ denotes the set of permutations $\pi$ such that if the labels of $\tilde{N}$ are permuted by $\pi$, then the resulting networks are the same.

We want to find $\vert F(\tilde{N})\vert$. Therefore, we need some notations.

First, a tree vertex $v$ is said to {\it root a subnetwork} if the set of all vertices $S$ which can be reached from $v$ (including $v$) has an induced subgraph $\tilde{N}(v)$ of $\tilde{N}$ which is connected to the set $V(\tilde{N})\setminus S$ only by the edge to $v$. $\tilde{N}(v)$ is called the {\it subnetwork} rooted at $v$. Moreover, we include the root into this definition which always roots a subnetwork, namely, $\tilde{N}$ itself.

Next,  a vertex $v$ which roots a subnetwork is called {\it symmetric} if the subnetwork $\tilde{N}(v)$ can be drawn in such a way that if it is reflected about the vertical line through $v$ by the angle $\pi$, then we obtain the same network if all labels of leaves are removed.

Likewise, we consider unordered pairs of tree vertices $\{v,w\}$ such that for the set of vertices $S$ which can be reached from $v$ and $w$ (including $v$ and $w$), the induced subnetwork of $\tilde{N}$ - denoted by $\tilde{N}(v,w)$ - is connected to $V(\tilde{N})\setminus S$ only via the edges to $v$ and $w$. (Note that $\tilde{N}(v,w)$ is not a phylogenetic networks since it has two roots.) Again, such a pair is said to be {\it symmetric} if $\tilde{N}(v,w)$ can be drawn such the line through $v$ and $w$ is a symmetry line.

Finally, we call a symmetric vertex $v$ (resp. symmetric pair of vertices $\{v,w\}$) {\it independent} if either (i) at least one symmetric vertex or at least one symmetric pair of $\tilde{N}(v)$ (resp. $\tilde{N}(v,w)$) does not lie on the symmetry line or (ii) at least one leaf which is not contained in a proper subnetwork of a symmetric vertex or symmetric pair of $\tilde{N}(v)$ (resp. $\tilde{N}(v,w)$) does not lie on the symmetry line; see Figure~\ref{fig-app} for examples. (Proper here means that the subnetwork is not equal to $\tilde{N}(v)$ (resp. $\tilde{N}(v,w)$).)

Denote now by $f$ the number of all independent vertices and pairs of vertices of $\tilde{N}$. Then,
\[
\vert F(\tilde{N})\vert=2^{f}.
\]
Thus, the answer to the above question is as follows.

\begin{theo}
Let $N$ be a phylogenetic network with $\ell$ leaves and $f$ the number of independent vertices and independent pairs of vertices of $N$. Then, the number of different leaf-labeled networks obtained from $N$ by labeling the leaves is $\ell!2^{-f}$.
\end{theo}

Moreover, this theorem has the following corollary which was used at the end of Section~\ref{normal}.

\begin{cor}\label{cor-app}
Let ${\mathcal C}$ be a class of leaf-labeled phylogenetic networks which is closed under permutations of the labels. Denote by ${\mathcal C}_\ell$ the networks from ${\mathcal C}$ with $\ell$ leaves. Then, $\vert{\mathcal C}_\ell\vert/\ell!$ is a fraction whose denominator is a power of $2$.
\end{cor}
\end{document}